\documentclass[a4paper, 10pt]{article}
\usepackage{graphicx}

\usepackage{alltt}

\usepackage{amsmath}
\usepackage{amssymb,amsmath,latexsym}
\usepackage[dvips]{color}
\usepackage[headings]{fullpage}
\usepackage{epic,epsfig,graphicx}
\usepackage{enumerate}
\usepackage{multirow}
\usepackage{setspace}
\usepackage{keyval}
\usepackage{float}

\usepackage{amssymb}

 \newtheorem{thm}{Theorem}[section]

 \newtheorem{lemma}{Lemma}[section]

 \newtheorem{defn}{Definition}[section]

 \newtheorem{rem}{Remark}[section]
  \newtheorem{assumption}{Assumption}[section]
 \newtheorem{note}{Note}[section]

 \numberwithin{equation}{section}

\begin{document}

\begin{center}
Nonlinear Analysis: Modelling and Control, Vol. vv, No. nn, YYYY\\
\copyright\ Vilnius University\\[24pt]
\LARGE
\textbf{Some asymptotic properties of SEIRS models with nonlinear incidence and random delays}\\[6pt]
\small
\textbf {Divine Wanduku\footnote{Corresponding author email: dwanduku@georgiasouthern.edu; wandukudivine@yahoo.com; Tel: +14073009605}, B. Oluyede\footnote{Email: boluyede@georgiasouthern.edu}}\\[6pt]
Department of Mathematical Sciences,
Georgia Southern University, \\65 Georgia Ave, Room 3042, Statesboro,
Georgia, 30460, U.S.A.\\[6pt]
Received: date\quad/\quad
Revised: date\quad/\quad
Published online: data
\end{center}

\begin{abstract}
This paper presents the dynamics of mosquitoes and humans, with general nonlinear incidence rate and multiple distributed delays for the disease. The model is a SEIRS system of delay differential equations. The normalized dimensionless version is derived; analytical techniques are applied to find conditions for deterministic extinction and permanence of disease. The BRN $R^{*}_{0}$ and ESPR $E(e^{-(\mu_{v}T_{1}+\mu T_{2})})$ are computed. Conditions for deterministic extinction and permanence are expressed in terms of $R^{*}_{0}$ and $E(e^{-(\mu_{v}T_{1}+\mu T_{2})})$, and applied to a P.vivax malaria scenario. Numerical results are given.
 \vskip 2mm


\end{abstract}

\nocite{2009ProcDETAp}
\begin{note}
\small{This arxiv paper published by Nonlinear Analysis: Modelling and Control  is the elaborate version with added biological insights which were removed to meet space restrictions of the journal. Thanks for reading. D.W.}
\end{note}
\section{Introduction\label{ch1.sec0}}

Malaria has exhibited an increasing alarming high mortality rate between 2015 and 2016. In fact, the latest WHO-\textit{World Malaria Report 2017} \cite{WHO-new} estimates a total of 216 million cases of malaria from 91 countries in 2016, which constitutes a 5 million increase in the total malaria cases from the malaria statistics obtained previously in 2015. Moreover, the total death count was 445000, and sub-Saharan Africa accounts for 90\% of the total estimated malaria cases. This rising trend in the malaria data, signals a need for more learning about the disease, improvement of the existing control strategies and equipment, and also a need for more advanced  resources etc. to fight  and eradicate, or ameliorate the burdens of malaria.

 Malaria and other mosquito-borne diseases such as dengue fever, yellow fever, zika fever, lymphatic filariasis etc. exhibit some  unique biological features. For instance, the incubation of the disease requires two hosts - the mosquito vector and human hosts,  which may be either  directly involved in a full life cycle of the infectious agent consisting of two separate and independent segments of sub-life cycles, which are completed separately inside the two hosts, or directly involved in two separate and independent half-life cycles of the infectious agent in the hosts. Therefore, there is a total latent time lapse of disease incubation which extends over the two segments of delay incubation times namely: (1) the incubation period of the infectious agent ( or the half-life cycle) inside the vector, and (2) the incubation period of the infectious agent (or the other half-life cycle) inside the human being (cf.\cite{WHO,CDC}). In fact,  the malaria plasmodium undergoes the first developmental half-life cycle called the \textit{sporogonic cycle} inside the female \textit{Anopheles} mosquito lasting approximately $10-18$ days, following a successful blood meal from an infectious human being through a mosquito bite. Moreover, the mosquito  becomes infectious.  The parasite  completes the second developmental half-life cycle called the \textit{exo-erythrocytic cycle} lasting about 7-30 days inside the exposed human being\cite{WHO,CDC}, whenever the parasite is transferred to  human being in the process of the infectious mosquito foraging for another blood meal.

The exposure and successful recovery from a malaria parasite, for example, \textit{falciparum vivae} induces natural immunity against the disease which can protect against subsequent severe outbreaks of the disease. Moreover, the effectiveness and duration of the naturally acquired immunity against malaria is determined by several factors such as the  species and the frequency of  exposure to the parasites (cf.\cite{CDC,denise}).

 Compartmental mathematical epidemic dynamic models  have been used to investigate the dynamics of several different types of infectious diseases including malaria\cite{pang,hyun}. In general, these models are classified as SIS, SIR, SIRS, SEIRS,  and  SEIR etc.\cite{qliu,wanduku-fundamental,Wanduku-2017,sen,cooke-driessche} epidemic dynamic models depending on the compartments of the disease classes directly involved in the general disease dynamics.
%
 Many compartmental mathematical models with delays have been studied \cite{Wanduku-2017,wanduku-delay,cooke-driessche}.

 Some important investigations in the study of population dynamic models expressed as  systems of differential equations are the permanence, and extinction of disease in the population, and also stability of the equilibria  over sufficiently long time. Several papers in the literature\cite{zhien,wanbiao,tend}  have addressed these topics. The extinction of disease seeks to find conditions that are sufficient for the disease related classes in the population such as, the exposed and infectious classes, to become extinct over sufficiently long time. The permanence of disease also answers the question about whether a significant number of people in the disease related classes will remain over sufficiently long time.  Disease eradication or persistence of disease in the steady state population seeks to find conditions sufficient for the equilibria to be stable asymptotically.
%

  The primary objectives of this paper include, to investigate (1) the extinction,  and (2) the permanence of disease in a family of SEIRS epidemic models.
  In other words, we find conditions that are sufficient for a disease such as malaria, to become extinct from the population over time, and also conditions that cause the disease to be permanent in the population over time.

The rest of this paper is presented as follows:- in Section~\ref{ch1.sec0.sec0}, the mosquito-human models are derived. In Section~\ref{ch1.sec1}, some model validation and preliminary results are presented. In Section~\ref{ch1.sec3a}, the results for the permanence of the disease are presented. Moreover, simulation results for the permanence of the disease in the population are presented in Section~\ref{ch1.sec4}. In Section~\ref{ch1.sec2b}, the results for the extinction of the disease are presented. Moreover, the numerical simulation results for the extinction of disease are presented in Section~\ref{ch1.sec4}.
 \section{Derivation of the mosquito-host dynamics}\label{ch1.sec0.sec0}
 The following assumptions are made to derive the epidemic model. Ideas from \cite{baretta-takeuchi1} will be used to derive the model for the mosquito-human dynamics.
 %

  (A) There are delays in the disease dynamics, and the delays represent the incubation period of the infectious agents (plasmodium or dengue fever virus etc.) in the vector $T_{1}$, and in the human host $T_{2}$. The third delay represents the natural immunity period $T_{3}$, where the delays are random variables with densities $f_{T_{1}}, t_{0}\leq T_{1}\leq h_{1}, h_{1}>0$, and $f_{T_{2}}, t_{0}\leq T_{2}\leq h_{2}, h_{2}>0$ and $f_{T_{3}}, t_{0}\leq T_{3}<\infty$ (cf. \cite{wanduku-biomath}).

  %
(B) The vector (e.g. mosquito) population consists of two main classes namely: the susceptible vectors $V_{s}$ and the infectious vectors $V_{i}$. Moreover, it is assumed that the total  vector population denoted $V_{0}$ is constant at any time, that is, $V_{s}(t)+V_{i}(t)=V_{0}, \forall t\geq t_{0}$, where $V_{0}>0$ is a positive constant. The susceptible vectors $V_{s}$ are infected by infectious humans $\hat{I}$, and after the incubation period $T_{1}$, the exposed vector becomes infectious $V_{i}$. Moreover, there is homogenous mixing between the vector-host populations.  Therefore, the birth rate and death rate of the vectors are equal, and denoted $\hat{\mu}_{v}$. It is assumed that the turnover of the vector population is very high, and the total number of vectors $V_{0}$ at any time $t$, is very large, and as a result, $\hat{\mu}_{v}$ is sufficiently large number. In addition, it is assumed that the total vectors $V_{0}$ is exceedingly larger than the total humans present at any time $t$, denoted $\hat{N}((t), t\geq t_{0}$. That is, $V_{0}>>\hat{N}((t), t\geq t_{0}$.

(C) The humans consists of susceptible $(\hat{S})$, Exposed $(\hat{E})$, Infectious $(\hat{I})$ and removed $(\hat{R})$ classes. The susceptibles are infected by the infectious vectors $V_{i}$, and become exposed (E). The infectious agent incubates for $T_{2}$ time units, and the exposed individuals become infectious $\hat{I}$. The infectious class recovers from the disease with temporary or sufficiently long natural immunity and become $(\hat{R})$. Therefore, the total population present at time $t$, $\hat{N}(t)=\hat{S}(t)+\hat{E}(t)+\hat{I}(t)+\hat{R}(t),\forall t\geq t_{0}$.

  Furthermore, it is  assumed that the interaction between the infectious vectors $V_{i}$ and  susceptible humans $\hat{S}$ exhibits nonlinear behavior, due to the overcrowding of the vectors as described in  (B), and resulting in psychological effects on the susceptible individuals which lead to change of behavior that limits the disease transmission rate, and consequently in a nonlinear character for the incidence rate characterized by the nonlinear incidence function $G$. $G$ satisfies the conditions of Assumption~\ref{ch1.sec0.assum1}.
  \begin{assumption}\label{ch1.sec0.assum1}
\begin{enumerate}
  \item [$A1$]$G(0)=0$; $A2$: $G(I)$ is strictly monotonic on $[0,\infty)$; $A3$: $G\in C^{2}([0,\infty), [0,\infty))$, and  $G''(I)<0$;
   $A4$. $\lim_{I\rightarrow \infty}G(I)=C, 0\leq C<\infty$; 
    $A5$: $G(I)\leq I, \forall I>0$; $A6$
     \begin{equation}\label{ch1.sec4.lemma1.eq1a}
     \left(\frac{G(x)}{x}-\frac{G(y)}{y}\right)\left(G(x)-G(y)\right)\leq 0, \forall x, y\geq0.
   \end{equation}
\end{enumerate}
\end{assumption}
These assumptions form an extension of the assumptions in \cite{qliu,wanduku-extinct,wanduku-biomath}. Some examples of incidence functions include $G(x)=\frac{x}{1+\theta x}, \theta>0$ etc.

(D) There is  constant birthrate of humans $\hat{B}$ in the population, and all births are susceptible individuals. It is also assumed that the natural deathrate of human beings in the population is $\hat{\mu}$ and individuals die additionally due to disease related causes at the rate $\hat{d}$. From a biological point of view, the average lifespan of vectors $\frac{1}{\hat{\mu}_{v}}$, is much less than the average lifespan of a human being in the absence of disease $\frac{1}{\hat{\mu}}$. It follows that assuming exponential lifetime for all individuals (both vector and host) in the population, then the survival probabilities over the time intervals of length  $T_{1}=s\in [t_{0},h_{1}]$, and $T_{2}=s\in [t_{0},h_{2}]$,  satisfy
\begin{equation}\label{ch1.sec0.eqn0.eq1}
e^{-\hat{\mu}_{v}T_{1}}<<e^{-\hat{\mu} T_{1}}\quad and\quad e^{-\hat{\mu}_{v}T_{1}-\hat{\mu} T_{2}}<<e^{-\hat{\mu}(T_{1} +T_{2})}.
\end{equation}
%
%

Applying similar ideas in \cite{baretta-takeuchi1}, the vector dynamics from (A)-(D) follows the system
\begin{eqnarray}
dV_{s}(t)&=&[-\Lambda e^{-\hat{\mu}_{v}T_{1}}\hat{I}(t-T_{1})V_{s}(t-T_{1})-\hat{\mu}_{v}V_{s}(t)+\hat{\mu}_{v}(V_{s}(t)+V_{i}(t))]dt,\label{ch1.sec0.eq0.eq3}\\
dV_{i}(t)&=&[\Lambda e^{-\hat{\mu}_{v}T_{1}}\hat{I}(t-T_{1})V_{s}(t-T_{1})-\hat{\mu}_{v}V_{i}(t)]dt,\label{ch1.sec0.eq0.eq4}\\
V_{0}&=&V_{s}(t)+V_{i}(t),\forall t\geq t_{0}, t_{0}\geq 0,\label{ch1.sec0.eq0.eq5}
\end{eqnarray}
where $\Lambda$ is the effective disease transmission rate from an infectious human being to a susceptible vector. Observe that the incidence rate of the disease into the vector population $\Lambda e^{-\hat{\mu}_{v}T_{1}}\hat{I}(t-T_{1})V_{s}(t-T_{1})$ represents new infectious vectors occurring at time $t$, which became exposed at earlier time $t-T_{1}$, and surviving natural death over the incubation period $T_{1}$, with survival probability rate $e^{-\hat{\mu}_{v}T_{1}}$, and are infectious at time $t$.  The detailed host population dynamics is derived as follows.

At time $t$, it follows from (C) that when susceptible humans $\hat{S}$ and infectious vectors $V_{i}$ interact with $\hat{\beta}$ effective contacts per vector, per unit time, then under the assumption of homogenous mixing, the incidence rate of the disease into the human population is given by the term $\hat{\beta}\hat{S}(t)V_{i}(t)$.
With the assumption of crowding effects of the vector population, it follows from (C) that the incidence rate of the disease can be written as
\begin{equation}\label{ch1.sec0.eqn0}
\hat{\beta}\hat{S}(t)G(V_{i}(t)),
\end{equation}
 where $G$ is the nonlinear incidence function satisfying the conditions in Assumption~\ref{ch1.sec0.assum1}.

It follows easily (cf.\cite{wanduku-biomath}) from the assumptions (A)-(D), and  (\ref{ch1.sec0.eqn0}) that for $T_{j},j=1,2,3$ fixed in the population, the dynamics of malaria in the human population is given by the system 
\begin{eqnarray}
d\hat{S}(t)&=&\left[ \hat{B}-\hat{\beta} \hat{S}(t)G(V_{i}(t)) - \hat{\mu} \hat{S}(t)+ \hat{\alpha} \hat{I}(t-T_{3})e^{-\hat{\mu} T_{3}} \right]dt,\nonumber\\
&&\label{ch1.sec0.eq3.eq1}\\
d\hat{E}(t)&=& \left[ \hat{\beta} \hat{S}(t)G(V_{i}(t)) - \hat{\mu} \hat{E}(t)\right.\nonumber\\
&&\left.-\hat{\beta} \hat{S}(t-T_{2}) e^{-\hat{\mu} T_{2}}G(V_{i}(t-T_{2})) \right]dt,\label{ch1.sec0.eq4.eq1}\\
&&\nonumber\\
d\hat{I}(t)&=& \left[\hat{\beta} \hat{S}(t-T_{2}) e^{-\hat{\mu} T_{2}}G(V_{i}(t-T_{2}))- (\hat{\mu} +\hat{d}+ \hat{\alpha}) \hat{I}(t) \right]dt,\nonumber\\
&&\label{ch1.sec0.eq5.eq1}\\
d\hat{R}(t)&=&\left[ \hat{\alpha} \hat{I}(t) - \hat{\mu} \hat{R}(t)- \hat{\alpha} \hat{I}(t-T_{3})e^{-\hat{\mu} T_{3}} \right]dt.\label{ch1.sec0.eq6}
\end{eqnarray}
Furthermore, the incidence function $G$ satisfies the conditions in Assumption~\ref{ch1.sec0.assum1}. And the initial conditions are given in the following:
\begin{eqnarray}
&&\left(\hat{S}(t),\hat{E}(t), \hat{I}(t), \hat{R}(t)\right)
=\left(\varphi_{1}(t),\varphi_{2}(t), \varphi_{3}(t),\varphi_{4}(t)\right), t\in (-T_{max},t_{0}],\nonumber\\
&&\varphi_{k}\in \mathcal{C}((-T_{max},t_{0}],\mathbb{R}_{+}),\forall k=1,2,3,4, \nonumber\\
&&\varphi_{k}(t_{0})>0,\forall k=1,2,3,4,\quad and\quad \max_{t_{0}\leq T_{1}\leq h_{1}, t_{0}\leq T_{2}\leq h_{2}, T_{3}\geq t_{0}}{(T_{1}+ T_{2}, T_{3})}=T_{max}\nonumber\\
 \label{ch1.sec0.eq6.eq1}
\end{eqnarray}
where $\mathcal{C}((-T_{max},t_{0}],\mathbb{R}_{+})$ is the space of continuous functions with  the supremum norm
\begin{equation}\label{ch1.sec0.eq6.eq2}
||\varphi||_{\infty}=\sup_{ t\leq t_{0}}{|\varphi(t)|}.
\end{equation}

It is shown in the following that the vector-host dynamics in (\ref{ch1.sec0.eq0.eq3})-(\ref{ch1.sec0.eq0.eq5}) and (\ref{ch1.sec0.eq3.eq1})-(\ref{ch1.sec0.eq6.eq1}) lead to the malaria model in \cite{wanduku-biomath}, which omits the dynamics of the vector population,  under the assumptions (A)-(D).

Firstly, observe that the system (\ref{ch1.sec0.eq3.eq1})-(\ref{ch1.sec0.eq6.eq1}) satisfies   [Theorem~3.1, \cite{wanduku-biomath}], and the total human population $\hat{N}(t)=\hat{S}(t)+\hat{E}(t)+\hat{I}(t)+\hat{R}(t),\forall t\geq t_{0}$ obtained from  system (\ref{ch1.sec0.eq3.eq1})-(\ref{ch1.sec0.eq6.eq1})  with initially condition that satisfies $N(t_{0})\leq \frac{\hat{B}}{\hat{\mu}}$, must satisfy
\begin{equation}\label{ch1.sec0.eq0.eq0}
  \limsup_{t\rightarrow \infty}{\hat{N}(t)}=\frac{\hat{B}}{\hat{\mu}}.
\end{equation}

Therefore, the assumption (B) above, interpreted as $\frac{\hat{N}(t)}{V_{0}}<<1, \forall t\geq t_{0}$  implies that
  \begin{equation}\label{ch1.sec0.eq0.eq1}
  \limsup_{t\rightarrow \infty}{\hat{N}(t)}=\frac{\hat{B}}{\hat{\mu}},\quad and\quad \frac{\left(\frac{\hat{B}}{\hat{\mu}}\right)}{V_{0}}<<1.
\end{equation}
Define
\begin{equation}\label{ch1.sec0.eq0.eq2}
  \epsilon=\frac{\left(\frac{\hat{B}}{\hat{\mu}}\right)}{V_{0}},
\end{equation}
then from  (\ref{ch1.sec0.eq0.eq1})-(\ref{ch1.sec0.eq0.eq2}), it follows that $\epsilon=\frac{\left(\frac{\hat{B}}{\hat{\mu}}\right)}{V_{0}}<<1$.

Employing similar reason in \cite{baretta-takeuchi1}, define two natural dimensionless time scales $\eta$ and $\varrho$ for the joint vector-host dynamics (\ref{ch1.sec0.eq0.eq3})-(\ref{ch1.sec0.eq0.eq5}) and (\ref{ch1.sec0.eq3.eq1})-(\ref{ch1.sec0.eq6.eq1}) in the following.
 \begin{eqnarray}
  \eta&=& \left(\frac{\hat{B}}{\hat{\mu}}\right)\Lambda t, \label{ch1.sec0.eq6.eq3}\\
  \varrho&=& V_{0}\Lambda t.\label{ch1.sec0.eq6.eq4}
\end{eqnarray}
  Note that since  the total vector population $V_{0}$ from (B) above is constant, that is, $V_{s}(t)+V_{i}(t)=V_{0}, \forall t  \geq t_{0}$, and from (\ref{ch1.sec0.eq0.eq0}) and  [Theorem~3.1, \cite{wanduku-biomath}] the total human $0<\hat{N}(t)\leq \frac{\hat{B}}{\hat{\mu}}, \forall t \geq t_{0}$, whenever $\hat{N}(t_{0})\leq \frac{\hat{B}}{\hat{\mu}}$, then the time scales $\eta$ and $\varrho$ arise naturally to rescale the total vector and maximum total human populations $V_{0}$ and $\left(\frac{\hat{B}}{\hat{\mu}}\right)$, respectively, at any time.
 The time scale $\varrho$ is  "fast", and  $\eta$ is  "slow" (cf. \cite{baretta-takeuchi1}).

Therefore, from above, let
\begin{equation}\label{ch1.sec0.eq6.eq5}
  \hat{V}_{i}(t)=\frac{V_{i}(t)}{V_{0}},\quad and \quad   \hat{V}_{s}(t)=\frac{V_{s}(t)}{V_{0}},
\end{equation}
be the dimensionless  vector variables, and
\begin{equation}\label{ch1.sec0.eq6.eq6}
  {S}(t)=\frac{\hat{S}(t)}{\left(\frac{\hat{B}}{\hat{\mu}}\right)},{I}(t)=\frac{\hat{I}(t)}{\left(\frac{\hat{B}}{\hat{\mu}}\right)}, {E}(t)=\frac{\hat{E}(t)}{\left(\frac{\hat{B}}{\hat{\mu}}\right)},   {R}(t)=\frac{\hat{R}(t)}{\left(\frac{\hat{B}}{\hat{\mu}}\right)}\quad and \quad {N}(t)=\frac{\hat{N}(t)}{\left(\frac{\hat{B}}{\hat{\mu}}\right)},
\end{equation}
 be the dimensionless human variables. And since $0<\hat{N}(t)\leq \frac{\hat{B}}{\hat{\mu}}, \forall t \geq t_{0}$, whenever $\hat{N}(t_{0})\leq \frac{\hat{B}}{\hat{\mu}}$, it follows from (\ref{ch1.sec0.eq6.eq6}) that
 \begin{equation}\label{ch1.sec0.eq6.eq6.eq1}
 0<{S}(t)+{E}(t)+{I}(t)+{R}(t)={N}(t)\leq 1, \forall t\geq t_{0}.
 \end{equation}

Applying (\ref{ch1.sec0.eq6.eq5})-(\ref{ch1.sec0.eq6.eq6}) to (\ref{ch1.sec0.eq0.eq3})-(\ref{ch1.sec0.eq0.eq5}) leads to the following
\begin{eqnarray}
d\hat{V}_{i}(t)&=& \epsilon\left[ e^{-\hat{\mu}_{v}T_{1}}{I}(t-T_{1})\hat{V}_{s}(t-T_{1})-\frac{\hat{\mu}_{v}}{\Lambda\left(\frac{\hat{B}}{\hat{\mu}}\right)}\hat{V}_{i}(t)\right]d\varrho,\label{ch1.sec0.eq6.eq7}\\
d\hat{V}_{s}(t)&=&-d\hat{V}_{i}(t),\label{ch1.sec0.eq6.eq8}\\
1&=&\hat{V}_{s}(t)+\hat{V}_{i}(t),\forall t\geq t_{0}, t_{0}\geq 0.\label{ch1.sec0.eq6.eq9}
\end{eqnarray}
Observe from (\ref{ch1.sec0.eq6.eq6.eq1})-(\ref{ch1.sec0.eq6.eq9}) that for nonnegative values for the vector variables $\hat{V}_{i}(t)\geq 0, \hat{V}_{s}(t)\geq 0, \forall t\geq t_{0}$, and positive values for the human variables ${S}(t), {E}(t), {I}(t), {R}(t)>0, \forall t\geq t_{0} $, it is follows that
\begin{equation}\label{ch1.sec0.eq6.eq10}
  -\epsilon\frac{\hat{\mu}_{v}}{\Lambda\left(\frac{\hat{B}}{\hat{\mu}}\right)}\leq \frac{d\hat{V}_{i}(t)}{d\varrho}\leq  \epsilon e^{-\hat{\mu}_{v}T_{1}}.
\end{equation}

Thus, on the time scale $\varrho$ which is "fast", it is easy to see from (\ref{ch1.sec0.eq6.eq7})-(\ref{ch1.sec0.eq6.eq10}), that under the assumption that $\epsilon$ from (\ref{ch1.sec0.eq0.eq2}) is infinitesimally small, that is $\epsilon\rightarrow 0$, then
\begin{equation}\label{ch1.sec0.eq6.eq11}
  \frac{d\hat{V}_{i}(t)}{d\varrho}=-\frac{d\hat{V}_{s}(t)}{d\varrho}=0,
\end{equation}
 which implies that the dynamics of $\hat{V}_{i}$ and $\hat{V}_{s}$ behaves as in steady state. And thus, it follows  from (\ref{ch1.sec0.eq6.eq7})-(\ref{ch1.sec0.eq6.eq11}) that
 \begin{eqnarray}
    \hat{V}_{i}(t)&=& \frac{e^{-\hat{\mu}_{v}T_{1}}}{\hat{\mu}_{v}}\Lambda \left(\frac{\hat{B}}{\hat{\mu}}\right){I}(t-T_{1})\hat{V}_{s}(t-T_{1}),\nonumber  \\
   1 &=& \hat{V}_{s}(t)+\hat{V}_{i}(t).\label{ch1.sec0.eq6.eq12}
 \end{eqnarray}
 It follows further from (\ref{ch1.sec0.eq6.eq12}) that
 \begin{equation}\label{ch1.sec0.eq6.eq13}
   \hat{V}_{s}(t)=\frac{1}{1+\frac{e^{-\hat{\mu}_{v}T_{1}}}{\hat{\mu}_{v}}\Lambda \left(\frac{\hat{B}}{\hat{\mu}}\right){I}(t-T_{1})\hat{V}_{s}(t-T_{1})}.
 \end{equation}
 For sufficiently large value of the birth-death rate $\hat{\mu}_{v}$ (see assumption (B)), such that $\hat{\mu}_{v}e^{\hat{\mu}_{v}T_{1}}>>\Lambda \left(\frac{\hat{B}}{\hat{\mu}}\right)$, then it follows from (\ref{ch1.sec0.eq6.eq13}) that $\hat{V}_{s}(t)\approx 1$, and consequently from (\ref{ch1.sec0.eq6.eq9}) and (\ref{ch1.sec0.eq6.eq5}), $V_{s}(t)\approx V_{0}$. Moreover, it follows further from (\ref{ch1.sec0.eq6.eq12}) that \begin{equation}\label{ch1.sec0.eq6.eq14}
  \hat{V}_{i}(t)\approx \frac{e^{-\hat{\mu}_{v}T_{1}}}{\hat{\mu}_{v}}\Lambda \left(\frac{\hat{B}}{\hat{\mu}}\right){I}(t-T_{1}),
\end{equation}
 and equivalently from (\ref{ch1.sec0.eq6.eq5})-(\ref{ch1.sec0.eq6.eq6}) that (\ref{ch1.sec0.eq6.eq14}) can be rewritten as follows
 \begin{equation}\label{ch1.sec0.eq6.eq15}
   V_{i}(t)\approx \frac{e^{-\hat{\mu}_{v}T_{1}}}{\hat{\mu}_{v}}\Lambda V_{0}\hat{I}(t-T_{1}).
 \end{equation}
 While on the fast scale $\varrho$ the term $\hat{I}(t-T_{1})$  behaves as the steady state, on the slow scale $\eta$, it is expected to still be evolving.  In the following, using (\ref{ch1.sec0.eq6.eq5})-(\ref{ch1.sec0.eq6.eq6}), the dynamics for the human population in (\ref{ch1.sec0.eq3.eq1})-(\ref{ch1.sec0.eq6.eq1}) is nondimensionalized with respect to the slow time scale $\eta$ in (\ref{ch1.sec0.eq6.eq3}).

  Without loss of generality(as it is usually the case e.g. $G(x)=\frac{x}{1+\alpha x}$, $G(x)=\frac{x}{1+\alpha x^{2}}$), it is assumed that on the $\eta$ timescale, the nonlinear term $G(V_{i}(t))$ expressed as $G(V_{0}\hat{V}_{i}(\eta))$, can be rewritten from (\ref{ch1.sec0.eq6.eq15})  as
 \begin{equation}\label{ch1.sec0.eq6.eq16}
   G(V_{0}\hat{V}_{i}(\eta))\equiv \frac{\Lambda V_{0}\left(\frac{\hat{B}}{\hat{\mu}}\right)}{\hat{\mu}_{v}}\hat{G}(\hat{V}_{i}(\eta))e^{-\hat{\mu}_{v}T_{1}},
 \end{equation}
by factoring a constant term $\frac{\Lambda V_{0}\left(\frac{\hat{B}}{\hat{\mu}}\right)}{\hat{\mu}_{v}}$, and the function $\hat{G}$ carries all the properties of Assumption~\ref{ch1.sec0.assum1}.
Thus, from the above and (\ref{ch1.sec0.eq6.eq15}), the  system (\ref{ch1.sec0.eq3.eq1})-(\ref{ch1.sec0.eq6.eq1}) is rewritten  in dimensionless form as follows:
\begin{eqnarray}
  d{S}(\eta) &=& [B-\beta {S}(\eta)\hat{G}({I}(\eta-T_{1\eta}))e^{-\mu_{v}T_{1\eta}}-\mu {S}(\eta)+\alpha I(\eta -T_{3\eta})e^{-\mu T_{3\eta}}]d\eta,\nonumber\\
  && \label{ch1.sec0.eq6.eq17}\\
   d{E}(\eta) &=& [\beta {S}(\eta)\hat{G}({I}(\eta-T_{1\eta}))e^{-\mu_{v}T_{1\eta}}-\mu {E}(\eta)\nonumber\\
   &&-\beta {S}(\eta-T_{2\eta})\hat{G}({I}(\eta-T_{1\eta}-T_{2\eta}))e^{-\mu_{v}T_{1\eta}-\mu T_{2\eta}}]d\eta, \label{ch1.sec0.eq6.eq18}\\
 d{I}(\eta) &=& [\beta {S}(\eta-T_{2\eta})\hat{G}({I}(\eta-T_{1\eta}-T_{2\eta}))e^{-\mu_{v}T_{1\eta}-\mu T_{2\eta}}-\mu {I}(\eta)\nonumber\\
 &&-(\mu+d+\alpha){I}(\eta)]d\eta,\label{ch1.sec0.eq6.eq19}\\
 d{R}(\eta)&=&[\alpha {I}(\eta)-\mu {R}(\eta)-\alpha I(\eta -T_{3\eta})e^{-\mu T_{3\eta}}]d\eta,\label{ch1.sec0.eq6.eq19}
\end{eqnarray}
where
\begin{eqnarray}
   && B=\frac{\hat{B}}{\left(\frac{\hat{B}}{\hat{\mu}}\right)^{2}\Lambda},\quad \beta=\frac{\hat{\beta}V_{0}}{\hat{\mu}_{v}},\quad \mu=\frac{\hat{\mu}}{\left(\frac{\hat{B}}{\hat{\mu}}\right)\Lambda},\quad \alpha=\frac{\hat{\alpha}}{\left(\frac{\hat{B}}{\hat{\mu}}\right)\Lambda}\nonumber\\
   &&\mu_{v}=\frac{\hat{\mu_{v}}}{\left(\frac{\hat{B}}{\hat{\mu}}\right)\Lambda},\quad d=\frac{\hat{d}}{\left(\frac{\hat{B}}{\hat{\mu}}\right)\Lambda},\quad  T_{j\eta}=\left(\frac{\hat{B}}{\hat{\mu}}\right)\Lambda T_{j},\forall j=1,2,3.
   \label{ch1.sec0.eq6.eq20}
\end{eqnarray}
The system (\ref{ch1.sec0.eq6.eq17})-(\ref{ch1.sec0.eq6.eq19}) describes the dynamics of malaria  on the slow scale $\eta$. Furthermore, moving forward, the analysis of the model (\ref{ch1.sec0.eq6.eq17})-(\ref{ch1.sec0.eq6.eq19}) is considered only on the $\eta$ timescale.  To reduce heavy notation, the following substitutions are made.
Substitute $t$ for $\eta$, and the delays $T_{j},\forall j=1,2,3$ will substitute $T_{j\eta}, \forall j=1,2,3$. Moreover, since the delays are are distributed with density functions $f_{T_{j}},\forall j=1,2,3$, it follows from (A)-(D), (\ref{ch1.sec0.eq6.eq17})-(\ref{ch1.sec0.eq6.eq19}) and (\ref{ch1.sec0.eq6.eq1}) that the expected SEIRS model for malaria is given as follows:
 \begin{eqnarray}
dS(t)&=&\left[ B-\beta S(t)\int^{h_{1}}_{t_{0}}f_{T_{1}}(s) e^{-\mu_{v} s}G(I(t-s))ds - \mu S(t)+ \alpha \int_{t_{0}}^{\infty}f_{T_{3}}(r)I(t-r)e^{-\mu r}dr \right]dt,\nonumber\\
&&\label{ch1.sec0.eq3}\\
dE(t)&=& \left[ \beta S(t)\int^{h_{1}}_{t_{0}}f_{T_{1}}(s) e^{-\mu_{v} s}G(I(t-s))ds - \mu E(t)\right.\nonumber\\
&&\left.-\beta \int_{t_{0}}^{h_{2}}f_{T_{2}}(u)S(t-u)\int^{h_{1}}_{t_{0}}f_{T_{1}}(s) e^{-\mu_{v} s-\mu u}G(I(t-s-u))dsdu \right]dt,\nonumber\\
&&\label{ch1.sec0.eq4}\\
dI(t)&=& \left[\beta \int_{t_{0}}^{h_{2}}f_{T_{2}}(u)S(t-u)\int^{h_{1}}_{t_{0}}f_{T_{1}}(s) e^{-\mu_{v} s-\mu u}G(I(t-s-u))dsdu- (\mu +d+ \alpha) I(t) \right]dt,\nonumber\\
&&\label{ch1.sec0.eq5}\\
dR(t)&=&\left[ \alpha I(t) - \mu R(t)- \alpha \int_{t_{0}}^{\infty}f_{T_{3}}(r)I(t-r)e^{-\mu s}dr \right]dt,\label{ch1.sec0.eq6}
\end{eqnarray}
where the initial conditions are given in the following: let $h= h_{1}+ h_{2}$ and define
\begin{eqnarray}
&&\left(S(t),E(t), I(t), R(t)\right)
=\left(\varphi_{1}(t),\varphi_{2}(t), \varphi_{3}(t),\varphi_{4}(t)\right), t\in (-\infty,t_{0}],\nonumber\\
&&\varphi_{k}\in U C_{g}\subset \mathcal{C}((-\infty,t_{0}],\mathbb{R}_{+}),\forall k=1,2,3,4,\quad \varphi_{k}(t_{0})>0,\forall k=1,2,3,4,\nonumber\\
 \label{ch1.sec0.eq06a}
\end{eqnarray}
where $UC_{g}$ is some fading memory sub Banach space of the Banach space  $\mathcal{C}((-\infty,t_{0}],\mathbb{R}_{+})$ endowed with the norm
\begin{equation}\label{ch1.sec0.eq06b}
\|\varphi\|_{g}=\sup _{t \leq t_{0}} \frac{|\varphi(t)|}{g(t)},
\end{equation}
and $g$ is some continuous function with the following properties: (P1.) $g\left(\left(-\infty, t_{0}\right]\right) \subseteq[1, \infty)$, non-increasing, and $g(t_{0})=1$; (P2.) $ \lim _{u \rightarrow t_{0}^{-}} \frac{g(t+u)}{g(t)}=1$, uniformly on $[t_{0},\infty)$; $ \lim _{t \rightarrow-\infty} g(t)=\infty$. An example of such a function is $g(t)=e^{-at}, a>0$ (cf. \cite{kuang}). Note that for any $g$ satisfying (P1.)-(P2.). the Banach space $\mathcal{C}((-\infty,t_{0}],\mathbb{R}_{+})$ is continuously embedded in $UC_{g}$ which allows structural properties for $\mathcal{C}((-\infty,t_{0}],\mathbb{R}_{+})$ with the uniform norm to hold in $UC_{g}$ with $||.||_{g}$ norm. Moreover, $\varphi \in U C_{g}, \exists g $ if and only if $||\varphi||_{g}<\infty $ and $\frac{|\varphi(t)|}{g(t)}$ is uniformly continuous on $(-\infty, t_{0}]$.
Also, the function $G$  in (\ref{ch1.sec0.eq3})-(\ref{ch1.sec0.eq6}) satisfies the conditions of Assumption~\ref{ch1.sec0.assum1}.

Observe  (\ref{ch1.sec0.eq3})-(\ref{ch1.sec0.eq6}) is similarly structured exactly as [(2.8)-(2.11), \cite{wanduku-biomath}]. Furthermore,
 the equations for $E$ and $R$ decouple from (\ref{ch1.sec0.eq3})-(\ref{ch1.sec0.eq6}). Therefore, the results are exhibited for the decoupled system (\ref{ch1.sec0.eq3}) and (\ref{ch1.sec0.eq5}) containing equations for $S$ and $I$. 
\begin{equation}\label{ch1.sec0.eq13b}
    Y(t)=(S(t), E(t), I(t), R(t))^{T},X(t)=(S(t),E(t),I(t))^{T},\quad\textrm{and}\quad N(t)=S(t)+ E(t)+ I(t)+ R(t).
\end{equation}

Whilst permanence or extinction has been investigated in  some delay type systems ( cf.\cite{zhien, wanbiao,tend}), the permanence  and extinction in the sense of \cite{zhien}  in systems with multiple random delays  is underdeveloped in the literature.
Furthermore,  as far as we know no other paper has addressed extinction and persistence of malaria in a mosquito-human population dynamics involving delay differential equations in the line of thinking of \cite{zhien,wanbiao}. We recall the following definition from \cite{zhien,xin}.
\begin{defn}\label{ch1.sec0.eq13b.def1}
\item[1.] A population $x(t)$ is called strongly permanent if
$\liminf_{t\rightarrow +\infty} x(t)>0$;
\item[2.] $x(t)$ is said to go extinct if
$\lim_{t\rightarrow +\infty}x(t)=0$.
\item[3.] $x(t)$ is said to be weakly permanent in the mean if
$\limsup_{t\rightarrow +\infty}\frac{1}{t}\int^{t}_{0}x(s)ds>0$.
\item[4.] $x(t)$ is said to be strongly permanent in the mean if
$\liminf_{t\rightarrow +\infty}\frac{1}{t}\int^{t}_{0}x(s)ds>0$.
\item[5.] $x(t)$ is said to  be stable in the mean if $\lim_{t\rightarrow \infty}{\frac{1}{t}\int^{t}_{t_{0}}x(s)}ds=c>0$.
\end{defn}
\section{Model validation results\label{ch1.sec1}}
The consistency results for the system (\ref{ch1.sec0.eq3})-(\ref{ch1.sec0.eq6}) are given. Some ideas from \cite{wanduku-biomath} using the dimensionless parameters (\ref{ch1.sec0.eq6.eq20}), are applied to the new model (\ref{ch1.sec0.eq3})-(\ref{ch1.sec0.eq6}). Observe from (\ref{ch1.sec0.eq6.eq20}) that expression $\frac{B}{\mu}$ simplifies to 1, and this is emphasized as $\frac{B}{\mu}\equiv 1$.
\begin{thm}\label{ch1.sec2b.lemma1.thm1}
For the given initial conditions (\ref{ch1.sec0.eq06a})-(\ref{ch1.sec0.eq06b}), the system (\ref{ch1.sec0.eq3})-(\ref{ch1.sec0.eq6}) has a unique positive solution $Y(t)\in \mathbb{R}_{+}^{4}$.
Moreover,
\begin{equation}\label{ch1.sec1.thm1a.eq0}
\limsup_{t\rightarrow \infty} N(t)\leq S^{*}_{0}=\frac{B}{\mu}\equiv 1.
   \end{equation}
  Furthermore, there is a positive self invariant space for the system denoted $D(\infty)=\bar{B}^{(-\infty, \infty)}_{\mathbb{R}^{4}_{+},}\left(0,\frac{B}{\mu}\equiv 1\right) $, where $D(\infty)$ is the closed unit ball in $\mathbb{R}^{4}_{+}$ centered at the origin with radius $\frac{B}{\mu}\equiv 1$ containing all positive solutions defined over $(-\infty,\infty )$.
  \end{thm}
  Proof:\\
  The proof of this result is standard and easy to follow applying the notations (\ref{ch1.sec0.eq13b}) to the system (\ref{ch1.sec0.eq3})-(\ref{ch1.sec0.eq6}).
\begin{thm}\label{ch1.sec2b.lemma1}
The feasible region for the unique positive solutions $Y(t), t\geq t_{0}$ of the  system (\ref{ch1.sec0.eq3})-(\ref{ch1.sec0.eq6}) in the phase plane that lie in the self-invariant unit ball $D(\infty)=\bar{B}^{(-\infty, \infty)}_{\mathbb{R}^{4}_{+},}\left(0,\frac{B}{\mu}\equiv 1\right) =\bar{B}^{(-\infty, \infty)}_{\mathbb{R}^{4}_{+},}\left(0, 1\right)$ for the system, also lie in a much smaller space $D^{expl}(\infty)\subset D(\infty)$, where
\begin{equation}\label{ch1.sec2b.lemma1.eq1}
  D^{expl}(\infty)=\left\{Y(t)\in \mathbb{R}^{4}_{+}:\frac{B}{\mu+d}\leq N(t)=S(t)+ E(t)+ I(t)+ R(t)\leq \frac{B}{\mu}, \forall t\in (-\infty, \infty) \right\}.
\end{equation}
 Moreover, the space $D^{expl}(\infty)$ is also self-invariant with respect to the system (\ref{ch1.sec0.eq3})-(\ref{ch1.sec0.eq6}).
\end{thm}
Proof:\\
Suppose $Y(t)\in D(\infty)$, then it follows from (\ref{ch1.sec0.eq3})-(\ref{ch1.sec0.eq6}) and (\ref{ch1.sec0.eq13b}) that the total population $N(t)=S(t)+E(t)+I(t)+R(t)$ satisfies the following inequality
\begin{equation}\label{ch1.sec2b.lemma1.proof.eq1}
[B-(\mu+d)N(t)]dt\leq dN(t)\leq [B-(\mu)N(t)]dt.
\end{equation}
It is easy to see from (\ref{ch1.sec2b.lemma1.proof.eq1}) that
\begin{equation}\label{ch1.sec2b.lemma1.proof.eq2}
\frac{B}{\mu+d}\leq \liminf_{t\rightarrow \infty}N(t)\leq \limsup_{t\rightarrow \infty}N(t)\leq \frac{B}{\mu},
\end{equation}
and (\ref{ch1.sec2b.lemma1.eq1}) follows immediately.
\begin{rem}
Theorem~\ref{ch1.sec2b.lemma1} signifies that every solution for (\ref{ch1.sec0.eq3})-(\ref{ch1.sec0.eq6}) that starts in the unit ball $D(\infty)$ in the phase plane, oscillates continuously inside $D(\infty)$. Moreover, if the solution oscillates and enters the space  $D^{expl}(\infty)=\bar{B}^{(-\infty, \infty)}_{\mathbb{R}^{4}_{+},}\left(0,\frac{B}{\mu}\equiv 1\right)\cap \left(\bar{B}^{(-\infty, \infty)}_{\mathbb{R}^{4}_{+},}\left(0,\frac{B}{\mu +d}\right)\right)^{c}$, the solution stays in $D^{expl}(\infty)$ for all time.

Biologically, observe that $\frac{B}{\mu }$ and $\frac{B}{\mu +d}$  represent the total births that occur over the average lifespans $\frac{1}{\mu}$ and $\frac{1}{\mu +d}$ of a human being in a malaria-free population and in a malaria- epidemic population, respectively. Thus, Theorem~\ref{ch1.sec2b.lemma1} signifies that when the population grows and enters a state for the total population $N(t)\in [\frac{B}{\mu+d },\frac{B}{\mu }]$, it stays within that range for all time. 
\end{rem}

Also, it is easy to see that the system (\ref{ch1.sec0.eq3})-(\ref{ch1.sec0.eq6}) has a DFE $E_{0}=(S^{*}_{0},0,0)=(\frac{B}{\mu}\equiv 1,0,0)=( 1,0,0)$. The basic reproduction number (BRN) for the disease when the delays in the system $T_{1}, T_{2}$ and $T_{3}$ are constant, is given  by
\begin{equation} \label{ch1.sec2.lemma2a.corrolary1.eq4}
\hat{R}^{*}_{0}=\frac{\beta }{(\mu+d+\alpha)}.
\end{equation}
Furthermore, when $\hat{R}^{*}_{0}<1$, then $E_{0}=X^{*}_{0}=(S^{*}_{0},0,0)=(1,0,0)$ is asymptotically stable, and  the disease can be eradicated from the population. Also, when the delays in the system  $T_{i}, i=1,2,3$  are random, and arbitrarily distributed, the BRN is proportional to 
\begin{equation}\label{ch1.sec2.theorem1.corollary1.eq3}
R_{0}\propto\frac{\beta }{(\mu+d+\alpha)}+\frac{\alpha}{(\mu+d+\alpha)},
\end{equation}
 And  malaria is eradicated from the system, whenever $R_{0}\leq 1$. 

The following result can be made about the nonzero steady state of the dimensionless system (\ref{ch1.sec0.eq3})-(\ref{ch1.sec0.eq6}), when Assumption~\ref{ch1.sec0.assum1} is satisfied.
\begin{thm}\label{ch1.sec2b.lemma1.thm2}
Let the conditions of Assumption~\ref{ch1.sec0.assum1} be satisfied. Suppose  $R_{0}>1$ (or $\hat{R}^{*}_{0}>1$) and the expected survival probability rate of the plasmodium satisfies
\begin{equation}
E(e^{-\mu_{v}T_{1}-\mu T_{2}})\geq \frac{R_{0}}{\left(R_{0}-\frac{\alpha}{\mu+d+\alpha}\right)G'(0)},
\end{equation}
  then there exists a nonzero endemic equilibrium $E_{1}=(S^{*}_{1}, E^{*}_{1}, I^{*}_{1})$ for the the dimensionless system (\ref{ch1.sec0.eq3})-(\ref{ch1.sec0.eq6}), where
  \begin{equation}\label{}
    E(e^{-\mu_{v}T_{1}-\mu T_{2}})=\int^{h_{2}}_{t_{0}}\int^{h_{1}}_{t_{0}}e^{-\mu_{v}s-\mu u}f_{T_{2}}(u)f_{T_{1}}(s)dsdu.
  \end{equation}
\end{thm}
Proof:\\
The dimensionless endemic equilibrium $E_{1}=(S^{*}_{1},  I^{*}_{1})$ of the decoupled  (\ref{ch1.sec0.eq3})-(\ref{ch1.sec0.eq6})  is a solution to the following system:
\begin{eqnarray}
 &&B-\beta E(e^{-\mu_{v} T_{1}}) SG(I)) - \mu S+ \alpha E(e^{-\mu T_{3}}) I =0,\label{ch1.sec3.thm1.proof.eq1}\\
&&\beta E(e^{-\mu_{v} T_{1}-\mu T_{2}}) SG(I)- (\mu +d+ \alpha) I=0.\label{ch1.sec3.thm1.proof.eq3}
\end{eqnarray}
Solving for $S$ from  (\ref{ch1.sec3.thm1.proof.eq3}) and substituting the result into (\ref{ch1.sec3.thm1.proof.eq1}), gives the following equation:
\begin{equation}\label{ch1.sec3.thm1.proof.eq3b}
H(I)=0
\end{equation}
 where,
\begin{equation}\label{ch1.sec3.thm1.proof.eq4}
H(I)= B-\frac{1}{E(e^{-\mu_{v}T_{1}+\mu T_{2}})}I\left[\frac{(\mu+d+\alpha)\mu}{\beta G(I)}+(\mu+d)E(e^{-\mu_{v} T_{1}})+\alpha E(e^{-\mu_{v} T_{1}})\left(1-E(e^{-\mu T_{2}})E(e^{-\mu T_{3}})\right)\right].
\end{equation}
 Note that $0<E(e^{-\mu T_{i}})\leq 1, i=1,2,3$, and $\lim_{I\rightarrow \infty} G(I)=C<\infty$,  hence for sufficiently large positive value of $I$,  $H(I)<0$. Furthermore, the derivative of $H(I)$ is given by
\begin{eqnarray}
H'(I)&=&-\frac{(\mu+d+\alpha)\mu}{\beta  E(e^{-\mu_{v} T_{1}-\mu T_{2}})} \frac{(G(I)-IG'(I))}{G^{2}(I)}\nonumber\\
&&-\frac{1}{E(e^{-\mu_{v} T_{1}-\mu T_{2}})}\left((\mu+d)E(e^{-\mu_{v} T_{1}})+\alpha E(e^{-\mu_{v} T_{1}})(1-E(e^{-\mu T_{2}})E(e^{-\mu T_{3}}))\right).\label{ch1.sec3.thm1.proof.eq5}
\end{eqnarray}
Assume without loss of generality that $G'(I)>0$. It follows from the other properties of $G$ in Assumption~\ref{ch1.sec0.assum1}, that is, $G(0)=0$, $G''(I)<0$,  that $(G(I)-IG'(I))>0$ and this further implies that $H'(I)<0$ for all $I> 0$. That is, $H(I)$ is a decreasing function  over all $I> 0$. Therefore, a positive root of the equation (\ref{ch1.sec3.thm1.proof.eq3b}) requires that $H(0)>0$. Observe from (\ref{ch1.sec3.thm1.proof.eq4}) and the dimensionless expressions in (\ref{ch1.sec0.eq6.eq20} ),
\begin{eqnarray}\label{ch1.sec3.thm1.proof.eq6}
&&H(0)=B\left(1- \frac{(\mu+d+\alpha)}{\beta  G'(0)E(e^{-\mu_{v} T_{1}-\mu T_{2}})}\right)
=B\left(1- \frac{1}{\left(R_{0}-\frac{\alpha}{\mu+d+\alpha}\right)G'(0)E(e^{-\mu_{v} T_{1}-\mu T_{2}})}\right)\nonumber\\
&&\geq B\left(1- \frac{1}{R_{0}}\right).
\end{eqnarray}
For $R_{0}>1$, it is easy to see that  $H(0)>0$.

%
The extinction of disease will be investigated in the neighborhood of the zero steady state $E_{0}$, and the permanence of disease will be investigated in the neighborhood of $E_{1}$.
\section{Extinction  of disease}\label{ch1.sec2b}
 In this section, the extinction of malaria from the system (\ref{ch1.sec0.eq3})-(\ref{ch1.sec0.eq6}) is investigated.
  Note, the decoupled system (\ref{ch1.sec0.eq3}) and (\ref{ch1.sec0.eq5}) is used.
The following lemma will be used to establish the extinction results.
\begin{lemma}\label{ch1.sec2b.lemma2}
Let the assumptions of Theorem~\ref{ch1.sec2b.lemma1} hold, and define the following Lyapunov functional in $D^{expl}(\infty)$,
\begin{eqnarray}
\tilde{V}(t)&=&V(t)+\beta\left[\int_{t_{0}}^{h_{2}}\int_{t_{0}}^{h_{1}}f_{T_{2}}(u)f_{T_{1}}(s)e^{-(\mu_{v}s+\mu u)}\int^{t}_{t-u}S(\theta)\frac{G(I(\theta-s))}{I(t)}d\theta dsdu\right.\nonumber\\
&&\left. +\int_{t_{0}}^{h_{2}}\int_{t_{0}}^{h_{1}}f_{T_{2}}(u)f_{T_{1}}(s)e^{-(\mu_{v}s+\mu u)}\int^{t}_{t-s}S(t)\frac{G(I(\theta))}{I(t)}d\theta ds du\right],
\label{ch1.sec2b.lemma2.eq1}
\end{eqnarray}
where $V(t)=\log{I(t)}$. It follows that
\begin{equation}\label{ch1.sec2b.lemma2.eq2}
\limsup_{t\rightarrow \infty}\frac{1}{t}\log{(I(t))}\leq \beta \frac{B}{\mu}E(e^{-(\mu_{v}T_{1}+\mu T_{2})})-(\mu+d+\alpha).
\end{equation}
\end{lemma}
Proof:\\
The differential operator $\dot{V}$ applied to the Lyapunov functional $\tilde{V}(t)$
with respect to the system (\ref{ch1.sec0.eq3}) leads to the following
\begin{equation}\label{ch1.sec2b.lemma2.proof.eq1}
  \dot{\tilde{V}}(t)=\beta\int_{t_{0}}^{h_{2}}f_{T_{2}}(u)\int^{h_{1}}_{t_{0}}f_{T_{1}}(s) e^{-(\mu_{v} s+\mu u)}S(t)\frac{G(I(t))}{I(t)}dsdu-(\mu+ d+ \alpha)
\end{equation}
Since $S(t), I(t)\in D^{expl}(\infty)$, and $G$ satisfies the conditions of Assumption~\ref{ch1.sec0.assum1}, it  follows easily from (\ref{ch1.sec2b.lemma2.proof.eq1}) that
\begin{equation}\label{ch1.sec2b.lemma2.proof.eq3}
\dot{\tilde{V}}(t)\leq \beta \frac{B}{\mu}E(e^{-(\mu_{v}T_{1}+\mu T_{2})})-(\mu+d+\alpha).
\end{equation}
Now, integrating both sides of  (\ref{ch1.sec2b.lemma2.proof.eq3}) over the interval $[t_{0},t]$, it follows from (\ref{ch1.sec2b.lemma2.proof.eq3})  and (\ref{ch1.sec2b.lemma2.eq1}) that
 \begin{eqnarray}
   \log{I(t)}&\leq&\tilde{V}(t)
    \leq\tilde{V}(t_{0})+\left[\beta \frac{B}{\mu}E(e^{-(\mu_{v}T_{1}+\mu T_{2})})-(\mu+d+\alpha)\right](t-t_{0}).\label{ch1.sec2b.lemma2.proof.eq4}
 \end{eqnarray}
 Diving both sides of (\ref{ch1.sec2b.lemma2.proof.eq4}) by $t$, and  taking the limit supremum  as $t\rightarrow\infty$, it is easy to see that (\ref{ch1.sec2b.lemma2.proof.eq4}) reduces to
 \begin{eqnarray}
   \limsup_{t\rightarrow\infty}\frac{1}{t}\log{I(t)}&\leq& \left[\beta \frac{B}{\mu}E(e^{-(\mu_{v}T_{1}+\mu T_{2})})-(\mu+d+\alpha)\right].\label{ch1.sec2b.lemma2.proof.eq5}
 \end{eqnarray}
 And the result (\ref{ch1.sec2b.lemma2.eq2}) follows immediately from  (\ref{ch1.sec2b.lemma2.proof.eq5}).

The extinction conditions for the infectious population over time are expressed in terms - (1) the BRN $R^{*}_{0}$ in (\ref{ch1.sec2.lemma2a.corrolary1.eq4}), and (2) the expected survival probability rate (ESPR) of the parasites $E(e^{-(\mu_{v}T_{1}+\mu T_{2})})$, also defined in [Theorem~5.1, \cite{wanduku-biomath}].
\begin{thm}\label{ch1.sec2b.thm1}
Suppose Lemma~\ref{ch1.sec2b.lemma2} is satisfied, and let the BRN $R^{*}_{0}$ be defined as in (\ref{ch1.sec2.lemma2a.corrolary1.eq4}). In addition, let one of the following conditions hold
\item[1.] $R^{*}_{0}\geq 1$ and $E(e^{-(\mu_{v}T_{1}+\mu T_{2})})<\frac{1}{R^{*}_{0}}$, or
\item[2.]$R^{*}_{0}<1$.\\
Then
\begin{equation}\label{ch1.sec2b.thm1.eq1}
\limsup_{t\rightarrow \infty}\frac{1}{t}\log{(I(t))}<-\lambda.
\end{equation}
where  $\lambda>0$ is some positive constant. In other words, $I(t)$ converges to zero exponentially.
\end{thm}
Proof:\\
Suppose Theorem~\ref{ch1.sec2b.thm1}~[1.] holds, then from (\ref{ch1.sec2b.lemma2.eq2}),
\begin{equation}\label{ch1.sec2b.thm1.eq1.proof.eq1}
\limsup_{t\rightarrow \infty}\frac{1}{t}\log{(I(t))}< \beta\frac{B}{\mu}\left(E(e^{-(\mu_{v}T_{1}+\mu T_{2})})- \frac{1}{R^{*}_{0}} \right)\equiv -\lambda,
\end{equation}
where the positive constant $\lambda>0$ is taken to be  as follows
\begin{equation}\label{ch1.sec2b.thm1.eq1.proof.eq1.eq1}
\lambda\equiv(\mu+d+\alpha)-\beta \frac{B}{\mu}E(e^{-(\mu_{v}T_{1}+\mu T_{2})})=\beta\frac{B}{\mu}\left( \frac{1}{R^{*}_{0}}-E(e^{-(\mu_{v}T_{1}+\mu T_{2})}) \right)>0.
\end{equation}
 Also, suppose Theorem~\ref{ch1.sec2b.thm1}~[2.] holds, then from (\ref{ch1.sec2b.lemma2.eq2}),
\begin{eqnarray}
\limsup_{t\rightarrow \infty}\frac{1}{t}\log{(I(t))}&\leq& \beta \frac{B}{\mu}E(e^{-(\mu_{v}T_{1}+\mu T_{2})})-(\mu+d+\alpha)\nonumber\\
&<& \beta \frac{B}{\mu}-(\mu+d+\alpha)= -(1-R^{*}_{0})(\mu+d+\alpha)\equiv -\lambda,\label{ch1.sec2b.thm1.eq1.proof.eq2}
\end{eqnarray}
where the positive constant $\lambda>0$ is taken to be  as follows
\begin{equation}\label{ch1.sec2b.thm1.eq1.proof.eq2.eq1}
\lambda\equiv(1-R^{*}_{0})(\mu+d+\alpha)>0.
\end{equation}
\begin{rem}
Theorem~\ref{ch1.sec2b.thm1},  and Theorem~\ref{ch1.sec2b.lemma1} signify that all trajectories of $(S(t), I(t))$  of the decoupled system (\ref{ch1.sec0.eq3}) and (\ref{ch1.sec0.eq5}), that  start in $D(\infty)$ and grow into $D^{expl}(\infty)\subset D(\infty)$ remain in $D^{expl}(\infty)$. Moreover, on the phase plane of $(S(t), I(t))$, the trajectory of the  infectious state $I(t), t\geq t_{0}$ ultimately turn to zero exponentially, whenever either the ESPR  $E(e^{-(\mu_{v}T_{1}+\mu T_{2})})<\frac{1}{R^{*}_{0}}$ for $R^{*}_{0}\geq 1$ , or whenever the BRN  $R^{*}_{0}<1$.  Furthermore, the  Lyapunov exponent (LE) from (\ref{ch1.sec2b.thm1.eq1}) is estimated by the  term $\lambda$, defined in (\ref{ch1.sec2b.thm1.eq1.proof.eq1.eq1}) and (\ref{ch1.sec2b.thm1.eq1.proof.eq2.eq1}).

  It follows from (\ref{ch1.sec2b.thm1.eq1}) that when either of the conditions in Theorem~\ref{ch1.sec2b.thm1}[1.-2.]  hold, then the $I(t)$ state dies out exponentially, whenever  $\lambda$ in (\ref{ch1.sec2b.thm1.eq1.proof.eq1.eq1}) and (\ref{ch1.sec2b.thm1.eq1.proof.eq2.eq1}) is positive, that is, $\lambda>0$. In addition, the rate of the exponential decay of each trajectories of $I(t)$ in each scenario of Theorem~\ref{ch1.sec2b.thm1}[1.-2.] is given by the  estimate $\lambda>0$ of the LE\footnote{Lyapunov exponent} in (\ref{ch1.sec2b.thm1.eq1.proof.eq1.eq1}) and (\ref{ch1.sec2b.thm1.eq1.proof.eq2.eq1}).

  The conditions in  Theorem~\ref{ch1.sec2b.thm1}[1.-2.] can also be interpreted as follows. Recall, the BRN   $R^{*}_{0}$ in (\ref{ch1.sec2.lemma2a.corrolary1.eq4}) (similarly in (\ref{ch1.sec2.theorem1.corollary1.eq3})) represents the expected number of secondary malaria cases that result from one infective placed in the disease free state $S^{*}_{0}=\frac{B}{\mu}\equiv 1$. Thus, $\frac{1}{R^{*}_{0}}=\frac{(\mu+d+\alpha)}{\beta S^{*}_{0}}$, for $R^{*}_{0}\geq 1$, represents the probability rate of  infectious persons in the secondary infectious population $\beta S^{*}_{0}$ leaving the infectious state, either through natural death $\mu$, diseases related death $d$, or recovery and acquiring natural immunity at the rate $\alpha$. Thus, $\frac{1}{R^{*}_{0}}$ is the effective probability rate of surviving infectiousness until recovery with acquisition of natural immunity. Moreover, $\frac{1}{R^{*}_{0}}$ is a probability measure provided $R^{*}_{0}\geq 1$.

  In addition, recall   Theorem\ref{ch1.sec2b.lemma1.thm2} asserts that when $R^{*}_{0}\geq 1$, and the ESPR $E(e^{-(\mu_{v}T_{1}+\mu T_{2})})$ is significantly large, then the outbreak of  malaria establishes a malaria endemic steady state population $E_{1}$. The conditions for extinction of disease in  Theorem~\ref{ch1.sec2b.thm1}[1.], that is $R^{*}_{0}\geq 1$ and $E(e^{-(\mu_{v}T_{1}+\mu T_{2})})<\frac{1}{R^{*}_{0}}$ suggest that in the event where $R^{*}_{0}\geq 1$, and the disease is aggressive, and likely to establish an endemic steady state population, if the expected survival probability rate $E(e^{-(\mu_{v}T_{1}+\mu T_{2})})$ of the malaria parasites over their complete life cycle of length $T_{1}+T_{2}$, is less than $\frac{1}{R^{*}_{0}}$- the effective probability rate of surviving infectiousness until recovery with natural immunity, then the malaria epidemic fails to establish an endemic steady state, and as a result, the disease ultimately dies out at an exponential rate $\lambda$ in (\ref{ch1.sec2b.thm1.eq1.proof.eq1.eq1}).
%

   In the event where  $R^{*}_{0}< 1$ in  Theorem~\ref{ch1.sec2b.thm1}[2.], extinction of disease occurs exponentially over sufficiently long time, regardless of the survival of the parasites. Moreover, the rate of extinction is  $\lambda$ in (\ref{ch1.sec2b.thm1.eq1.proof.eq2.eq1}).
\end{rem}
\section{Persistence of susceptibility and stability of zero equilibrium }
Theorem~\ref{ch1.sec2b.thm1} characterizes the behavior of the trajectories of the $I(t)$ coordinate of the solution $(S(t),I(t))$ of the decoupled system (\ref{ch1.sec0.eq3}) and (\ref{ch1.sec0.eq5}) in the phase plane. The question remains about how the trajectories for the $S(t)$ behave asymptotically in the phase plane.
   %

 Using Definition~\ref{ch1.sec0.eq13b.def1}[3-5],  the following result describes the average behavior of the trajectories of  $S(t)$ over sufficiently long time, and also states conditions for the  stability of the disease-free equilibrium (DFE) of the decoupled system  $E_{0}=(S^{*}_{0},0)=(1,0)$, whenever  Theorem~\ref{ch1.sec2b.thm1} holds.
\begin{thm}\label{ch1.sec2b.thm2}
Suppose any of the conditions in the hypothesis of Theorem~\ref{ch1.sec2b.thm1}[1.-2.] are satisfied. It follows that in $D^{expl}(\infty)$, the trajectories of the susceptible state $S(t)$ of the decoupled system (\ref{ch1.sec0.eq3}) and (\ref{ch1.sec0.eq5}), satisfy
\begin{equation}\label{ch1.sec2b.thm2.eq1}
\lim_{t\rightarrow \infty}\frac{1}{t}\int_{t_{0}}^{t}S(\xi)d\xi=\frac{B}{\mu}\equiv 1.
\end{equation}
That is, the susceptible state is strongly persistent over long-time in the mean (see definition of  persistence in the mean Definition~\ref{ch1.sec0.eq13b.def1}[3-4]). Moreover, it is  stable in the mean, and  the average value of the susceptible state over sufficiently long time is equal to $S(t)=S^{*}_{0}=\frac{B}{\mu}$, obtained when the system is in steady state.
\end{thm}
Proof:\\
Suppose either of the conditions in Theorem~\ref{ch1.sec2b.thm1}[1.-2.] hold, then  it follows clearly from Theorem~\ref{ch1.sec2b.thm1} that  for every $\epsilon>0$, there is a positive constant $K_{1}(\epsilon)\equiv K_{1}>0$, such that
\begin{equation}\label{ch1.sec2b.thm2.proof.eq2}
  I(t)<\epsilon,\quad  \textrm{whenever $t>K_{1}$}.
\end{equation}
It follows from (\ref{ch1.sec2b.thm2.proof.eq2}) that
\begin{equation}\label{ch1.sec2b.thm2.proof.eq3}
  I(t-s)<\epsilon,\quad \textrm{whenever $t>K_{1}+h_{1},\forall s\in [t_{0}, h_{1}]$}.
\end{equation}
In $D^{expl}(\infty)$, define
\begin{equation}\label{ch1.sec2b.thm2.proof.eq4}
V_{1}(t)=S(t)+\alpha \int_{t_{0}}^{\infty} f_{T_{3}}(r)e^{\mu r}\int_{t-r}^{t}I(\theta)d\theta dr.
\end{equation}
The differential operator $\dot{V}_{1}$ applied to the Lyapunov functional $V_{1}(t)$ in (\ref{ch1.sec2b.thm2.proof.eq4}) leads to the following
\begin{equation}\label{ch1.sec2b.thm2.proof.eq5}
  \dot{V}_{1}(t)=g(S, I)-\mu S(t),
\end{equation}
where
\begin{equation}\label{ch1.sec2b.thm2.proof.eq6}
  g(S,I)=B-\beta S(t)\int^{h_{1}}_{t_{0}}f_{T_{1}}(s) e^{-\mu_{v} s}G(I(t-s))ds+\alpha E(e^{-\mu T_{3}}) I(t).
\end{equation}
Estimating the right-hand-side of (\ref{ch1.sec2b.thm2.proof.eq5}) in $D^{expl}(\infty)$, and integrating over $[t_{0},t]$, it follows from (\ref{ch1.sec2b.thm2.proof.eq2})-(\ref{ch1.sec2b.thm2.proof.eq3}) that
\begin{eqnarray}
  V_{1}(t)&\leq& V_{1}(t_{0})+B(t-t_{0})+\int_{t_{0}}^{K_{1}}\alpha I(\xi)d\xi +\int_{K_{1}}^{t}\alpha I(\xi)d\xi -\mu\int_{t_{0}}^{t} S(\xi)d\xi,\nonumber\\
  &\leq& V_{1}(t_{0})+B(t-t_{0})+\alpha \frac{B}{\mu}(K_{1}-t_{0})  +\alpha (t-K_{1})\epsilon -\mu\int_{t_{0}}^{t} S(\xi)d\xi.\label{ch1.sec2b.thm2.proof.eq7}
\end{eqnarray}
Thus, dividing both sides of (\ref{ch1.sec2b.thm2.proof.eq7}) by $t$ and taking the limit supremum as $t\rightarrow \infty$, it follows that
\begin{equation}\label{ch1.sec2b.thm2.proof.eq10}
\limsup_{t\rightarrow \infty}\frac{1}{t}\int_{t_{0}}^{t}S(\xi)d\xi\leq \frac{B}{\mu}+ \frac{\alpha}{\mu}\epsilon.
\end{equation}
On the other hand, estimating $g(S,I)$ in (\ref{ch1.sec2b.thm2.proof.eq6}) from below and using the conditions of Assumption~\ref{ch1.sec0.assum1} and (\ref{ch1.sec2b.thm2.proof.eq3}), it is easy to see that  in $D^{expl}(\infty)$,
\begin{eqnarray}
  g(S,I) &\geq & B-\beta S(t)\int^{h_{1}}_{t_{0}}f_{T_{1}}(s) e^{-\mu_{v} s}(I(t-s))ds
   \geq B-\beta \frac{B}{\mu}E(e^{-\mu_{v}T_{1}})\epsilon,\forall t>K_{1}+h_{1},\nonumber\\
   &\geq&B-\beta \frac{B}{\mu}\epsilon.\label{ch1.sec2b.thm2.proof.eq11}
\end{eqnarray}
Moreover, for $t\in[t_{0}, K_{1}+h_{1}]$, then
\begin{equation}\label{ch1.sec2b.thm2.proof.eq11.eq1}
  g(S,I)\geq B-\beta \left(\frac{B}{\mu}\right)^{2}.
\end{equation}
Therefore, applying (\ref{ch1.sec2b.thm2.proof.eq11})-(\ref{ch1.sec2b.thm2.proof.eq11.eq1}) into (\ref{ch1.sec2b.thm2.proof.eq5}), then integrating both sides of (\ref{ch1.sec2b.thm2.proof.eq5}) over $[t_{0},t]$, and diving the result by $t$,  it is easy to see from (\ref{ch1.sec2b.thm2.proof.eq5}) that
\begin{equation}\label{ch1.sec2b.thm2.proof.eq12}
  \frac{1}{t}V_{1}(t)\geq  \frac{1}{t}V_{1}(t_{0})+B(1-\frac{t_{0}}{t}) -\frac{1}{t}\beta\left( \frac{B}{\mu}\right)^{2}(K_{1}+h_{1}-t_{0})-\beta\frac{B}{\mu}\epsilon[1-\frac{K_{1}+h_{1}}{t}] -\frac{1}{t}\mu\int_{t_{0}}^{t}S(\xi)d\xi.
\end{equation}
Observe that in  $D^{expl}(\infty)$, $\lim_{t\rightarrow \infty}\frac{1}{t}V_{1}(t)=0$, and $\lim_{t\rightarrow \infty}\frac{1}{t}V_{1}(t_{0})=0$.  Therefore, rearranging (\ref{ch1.sec2b.thm2.proof.eq12}), and taking the limit infinimum of both sides as $t\rightarrow \infty$, it is easy to see that
\begin{equation}\label{ch1.sec2b.thm2.proof.eq13}
 \liminf_{t\rightarrow \infty} \frac{1}{t}\int_{t_{0}}^{t}S(\xi)d\xi\geq   \frac{B}{\mu}-\frac{1}{\mu}\beta \frac{B}{\mu}\epsilon.
\end{equation}
It follows from (\ref{ch1.sec2b.thm2.proof.eq10}) and (\ref{ch1.sec2b.thm2.proof.eq13}) that
\begin{equation}\label{ch1.sec2b.thm2.proof.eq14}
\frac{B}{\mu}-\frac{1}{\mu}\beta \frac{B}{\mu}\epsilon\leq \liminf_{t\rightarrow \infty} \frac{1}{t}\int_{t_{0}}^{t}S(\xi)d\xi\leq \limsup_{t\rightarrow \infty}\frac{1}{t}\int_{t_{0}}^{t}S(\xi)d\xi\leq \frac{B}{\mu}+ \frac{\alpha}{\mu}\epsilon.
\end{equation}
Hence, for $\epsilon$ arbitrarily small,  the result in (\ref{ch1.sec2b.thm2.eq1}) follows immediately from (\ref{ch1.sec2b.thm2.proof.eq14}).
\begin{rem}\label{ch1.sec2b.thm2.proof.rem1a}
Theorem~\ref{ch1.sec2b.thm2} signifies that
 the DFE $E_{0}$ is strongly persistent and stable in the mean by Definition~\ref{ch1.sec0.eq13b.def1}[3-5]. That is, over sufficiently long time, on average the human population will be in the DFE $E_{0}$.
 Thus, the conditions in Theorem~\ref{ch1.sec2b.thm1} are sufficient for malaria to be eradicated from the population, when the population is in a steady state.
\end{rem}
The next, result confirms that not only is the zero equilibrium state of the decoupled system (\ref{ch1.sec0.eq3}) and (\ref{ch1.sec0.eq5})  $E_{0}=(S^{*}_{0},0)=(1,0)$ stable  and persistent on average over time, but also stable in the sense of Lyapunov.
\begin{thm}\label{ch1.sec2b.thm3}
Suppose any of the conditions in the hypothesis of Theorem~\ref{ch1.sec2b.thm1}[1.-2.] are satisfied. Also, suppose the conditions of Theorem~\ref{ch1.sec2b.thm2} hold. It follows that in $D^{expl}(\infty)$, the DFE $E_{0}=(S^{*}_{0},0)=(\frac{B}{\mu},0)=(1,0)$  is stable in the sense of Lyapunov.
\end{thm}
Proof:\\
It is left to show that every trajectory that starts near $E_{0}$ remains near $E_{0}$ asymptotically. Indeed, if the hypothesis of Theorem~\ref{ch1.sec2b.thm1}[1.-2.] holds, then all trajectories in the phase-plane for the infectious state $I(t)$ converge asymptotically and exponentially to $I^{*}_{0}=0$. It is left to show that if the trajectories of the susceptible state $S(t)$ from Theorem~\ref{ch1.sec2b.thm2} (\ref{ch1.sec2b.thm2.eq1}), converge asymptotically in the mean to $S^{*}_{0}=\frac{B}{\mu}$, then they  must remain asymptotically near $S^{*}_{0}=\frac{B}{\mu}$.

Indeed, if on the contrary, there exist a trajectory for $S(t)$ starting near $S^{*}_{0}=\frac{B}{\mu}=\equiv 1$ that does not stay near $S^{*}_{0}=\frac{B}{\mu}$ asymptotically, that is, suppose there exists some $\epsilon_{0}>0$ and $\delta (t_{0},\epsilon_{0})>0$, such that $||S(t_{0})-S^{*}_{0}||<\delta$, but $||S(t)-S^{*}_{0}||\geq \epsilon_{0}, \forall t\geq t_{0}$, then clearly from (\ref{ch1.sec2b.thm2.eq1}), either
\begin{equation}\label{ch1.sec2b.thm2.rem1.eq1}
 S^{*}_{0}=\lim_{t\rightarrow \infty}\frac{1}{t}\int_{t_{0}}^{t}S(\xi)d\xi\geq S^{*}_{0}+\epsilon_{0}\quad or\quad S^{*}_{0}=\lim_{t\rightarrow \infty}\frac{1}{t}\int_{t_{0}}^{t}S(\xi)d\xi\leq S^{*}_{0}-\epsilon_{0}.
\end{equation}
Thus, $\epsilon_{0}$ must be zero, otherwise  (\ref{ch1.sec2b.thm2.rem1.eq1}) is a contradiction. Hence, $E_{0}=(S^{*}_{0},0)$ is stable in the sense of Lyapunov.
\begin{rem}\label{ch1.sec2b.thm2.rem1}
Theorem~\ref{ch1.sec2b.thm2}, Theorem~\ref{ch1.sec2b.thm1}, and Theorem~\ref{ch1.sec2b.lemma1} signify that  all trajectories of  $(S(t), I(t)),t\geq t_{0}$ of the decoupled system (\ref{ch1.sec0.eq3}) and (\ref{ch1.sec0.eq5}) that start  in $D^{expl}(\infty)\subset D(\infty)$ remain bounded  in $D^{expl}(\infty)$. Moreover, the trajectories of  $I(t), t\geq t_{0}$ of the solution $(S(t), I(t)),t\geq t_{0}$ in phase plane, ultimately turn to zero exponentially, while trajectories of the susceptible state $S(t)$ persist strongly, and converge in the mean to the DFE $S^{*}_{0}=\frac{B}{\mu}$, whenever $E(e^{-(\mu_{v}T_{1}+\mu T_{2})})<\frac{1}{R^{*}_{0}}$, for $R^{*}_{0}\geq 1$, or whenever the basic production number satisfy $R^{*}_{0}<1$.
 %

Moreover, from Theorem~\ref{ch1.sec2b.thm3}, the DFE  $(S(t), I(t))=E_{0}=(S^{*}_{0},0)=(\frac{B}{\mu},0)=(1,0)$ is uniformly globally asymptotically stable.
Thus, the conditions in Theorem~\ref{ch1.sec2b.thm1} are strong disease eradication conditions.

The above observations suggest that over sufficiently long time, the population that remains will be all susceptible malaria-free people, and the population size will be averagely equal to the DFE $S^{*}_{0}=\frac{B}{\mu}$ of  (\ref{ch1.sec0.eq3}) and (\ref{ch1.sec0.eq5}).
\end{rem}
 \section{Permanence of infectivity near nonzero equilibrium\label{ch1.sec3a}}
 As remarked in Theorem~\ref{ch1.sec2b.lemma1.thm2}, when $R^{*}_{0}$ in (\ref{ch1.sec2.lemma2a.corrolary1.eq4}) satisfies $R^{*}_{0}>1$, the endemic equilibrium of the decoupled system (\ref{ch1.sec0.eq3}) and (\ref{ch1.sec0.eq5}) exists and  is denoted $E_{1}=(S^{*}_{1}, I^{*}_{1})$. In this section, conditions for $I(t)$ to be strongly persistent (Definition~\ref{ch1.sec0.eq13b.def1}[1]) in the neighborhood of $E_{1}$ are given.
 \begin{lemma}\label{ch1.sec4.lemma1}
   Suppose the conditions of Theorem~\ref{ch1.sec2b.lemma1} and Theorem~\ref{ch1.sec2b.lemma1.thm2}
    are satisfied, and let the nonlinear incidence function $G$ satisfy the assumptions of Assumption~\ref{ch1.sec0.assum1}.

   Then every positive solution $(S(t), I(t))\in D(\infty)$ of  the  decoupled  system (\ref{ch1.sec0.eq3}) and (\ref{ch1.sec0.eq5}) with initial conditions (\ref{ch1.sec0.eq06a}) and (\ref{ch1.sec0.eq06b}) satisfies the following conditions:
    \begin{equation}\label{ch1.sec4.lemma1.eq1}
    \liminf_{t\rightarrow \infty}{S(t)}\geq v_{1}\equiv \frac{B}{\mu+\beta G(S_{0})}\quad and\quad  \liminf_{t\rightarrow \infty}{I(t)}\geq v_{2}\equiv qI^{*}_{1}e^{-(\mu+d+\alpha)(\rho+1)h},
  \end{equation}
  where $h=h_{1}+h_{2}$,  and  $\rho>0$ is a suitable positive constant,  $S^{*}_{1}<\min\{S_{0}, S^{\vartriangle}\}$ and $0<q<\bar{q}<1$, given that,
  \begin{equation}\label{ch1.sec4.lemma1.eq2}
  \bar{q}=\frac{B\beta E(e^{-\mu T_{1}})G(I^{*}_{1})-\mu \alpha E(e^{-\mu T_{3}})I^{*}_{1}}{\left(B+\alpha E(e^{-\mu T_{3}})I^{*}_{1}\right)\beta I^{*}_{1}},\quad S^{\vartriangle}=\frac{B}{k}\left(1-e^{-k\rho h}\right), k=\mu +\beta G(q I^{*}_{1}).
  \end{equation}
 \end{lemma}
 Proof:\\
 Recall (\ref{ch1.sec1.thm1a.eq0}) asserts that for $N(t)=S(t)+E(t)+ I(t)+R(t)$, $\limsup_{t\rightarrow \infty} {N(t)}\leq S^{*}_{0}=\frac{B}{\mu}$. This implies that  $\limsup_{t\rightarrow \infty} {S(t)}\leq S^{*}_{0}\equiv 1$. This further implies that for any arbitrarily small $\epsilon>0$, there exists a sufficiently large $\Lambda>0$, such that
 \begin{equation}\label{ch1.sec4.lemma1.proof.eq1}
 I(t)\leq S^{*}_{0} +\varepsilon,\quad whenever, \quad t\geq \Lambda.
 \end{equation}
 Without loss of generality, let $\Lambda_{1}>0$ be sufficiently large such that
 \[t\geq \Lambda\geq \max_{(s, r)\in [t_{0}, h_{1}]\times [t_{0}, \infty)}{(\Lambda_{1}+s, \Lambda_{1}+r)}.\]
   It follows from Assumption~\ref{ch1.sec0.assum1}, (\ref{ch1.sec0.eqn0.eq1})  and (\ref{ch1.sec0.eq3}) that
 \begin{eqnarray}
 \frac{dS(t)}{dt}\geq B-\beta S(t)\int_{t_{0}}^{h_{1}}f_{T_{1}}(s)e^{-\mu s}G(S^{*}_{0}+\epsilon)ds-\mu S(t)
  \geq B-\left[\mu+ \beta G(S^{*}_{0}+\epsilon)\right]S(t).\label{ch1.sec4.lemma1.proof.eq2}
 \end{eqnarray}
 From (\ref{ch1.sec4.lemma1.proof.eq2}) it follows that
 \begin{equation}\label{ch1.sec4.lemma1.proof.eq3}
 S(t)\geq  \frac{B}{k_{1}} -\frac{B}{k_{1}}e^{-k_{1}(t-t_{0})}+S(t_{0}) e^{-k_{1}(t-t_{0})},
 \end{equation}
 where $k_{1}=\mu+ \beta G(S^{*}_{0}+\epsilon)$.

 It is easy to see from (\ref{ch1.sec4.lemma1.proof.eq3})
 \begin{equation}\label{ch1.sec4.lemma1.proof.eq4}
   \liminf_{t\rightarrow \infty}{S(t)}\geq  \frac{B}{\mu+\beta G(S_{0}+\epsilon)}.
 \end{equation}
 Since $\epsilon>0$ is arbitrarily small, then the first part of (\ref{ch1.sec4.lemma1.eq1}) follows immediately.

 In the following it is shown that $ \liminf_{t\rightarrow \infty}{I(t)}\geq v_{2}$. In order to establish this result, it is first proved that it is impossible that $I(t)\leq q I^{*}_{1}$ for sufficiently large  $t\geq t_{0}$, where $q\in(0, 1)$ is defined in the hypothesis. Suppose on the contrary there exists some sufficiently large $\Lambda_{0}> t_{0}>0$, such that $I(t)\leq q I^{*}_{1}, \forall t\geq \Lambda_{0}$. It follows from (\ref{ch1.sec0.eq3}) that
 \begin{eqnarray}
   S^{*}_{1} = \frac{B+\alpha E(e^{-\mu T_{3}})I^{*}_{1}}{\mu+ \beta E(e^{-\mu T_{1}})G(I^{*}_{1})} 
    = \frac{B}{\mu+\frac{B\beta E(e^{-\mu T_{1}})G(I^{*}_{1})-\mu \alpha E(e^{-\mu T_{3}})I^{*}_{1} }{B+\alpha E(e^{-\mu T_{3}})I^{*}_{1}}}.\label{ch1.sec4.lemma1.proof.eq5}
 \end{eqnarray}
 But, it can be easily seen from (\ref{ch1.sec0.eq3}) and (\ref{ch1.sec0.eq5}) that
 \begin{eqnarray}
    &&B\beta E(e^{-\mu T_{1}})G(I^{*}_{1})-\mu \alpha E(e^{-\mu T_{3}})I^{*}_{1}=\frac{\mu(\mu+d+\alpha)\left[S^{*}_{0}-\frac{\alpha E(e^{-\mu T_{3}})E(e^{-\mu T_{2}})}{(\mu +d+\alpha)}S^{*}_{1}\right]}{E(e^{-\mu T_{2}})S^{*}_{1}}I^{*}_{1}\nonumber \\
   &&\geq \frac{\mu(\mu+d+\alpha)(S^{*}_{0}-S^{*}_{1})}{E(e^{-\mu T_{2}})S^{*}_{1}}>0,\quad since\quad S^{*}_{0}=\frac{B}{\mu}\geq S^{*}_{1}.\label{ch1.sec4.lemma1.proof.eq6}
    \end{eqnarray}
   Therefore, from (\ref{ch1.sec4.lemma1.proof.eq5}), it follows that
   \begin{equation}\label{ch1.sec4.lemma1.proof.eq6}
     S^{*}_{1}<\frac{B}{\mu +\beta I^{*}_{1}q}\leq \frac{B}{\mu +\beta G(qI^{*}_{1})},
   \end{equation}
   where $0<q<\bar{q}$, and $\bar{q}$ is defined in (\ref{ch1.sec4.lemma1.eq2}).

   For all  vector values $ (s, r)\in [t_{0}, h_{1}]\times [t_{0}, \infty) $ define
   \begin{equation}
   \Lambda_{0,max}=\max_{(s, r)\in [t_{0}, h_{1}]\times [t_{0}, \infty)}{(\Lambda_{0}+s, \Lambda_{0}+r)},
   \end{equation}
It follows from Assumption~\ref{ch1.sec0.assum1}  and (\ref{ch1.sec0.eq3}) that for all $t\geq \Lambda_{0,max}$,
\begin{equation}\label{ch1.sec4.lemma1.proof.eq7}
   S(t)\geq  \frac{B}{k} -\frac{B}{k}e^{-k(t-\Lambda_{0,max})}+S(\Lambda_{0,max}) e^{-k(t-\Lambda_{0,max})},
\end{equation}
where $k$ is defined in (\ref{ch1.sec4.lemma1.eq2}).
 For $t\geq \Lambda_{0,max}+ \rho h$,  where $ h=h_{1}+h_{2}$, and $\rho>0$ is sufficiently large,  it follows from (\ref{ch1.sec4.lemma1.proof.eq7}) that
 \begin{equation}\label{ch1.sec4.lemma1.proof.eq8}
   S(t)\geq \frac{B}{k}\left[1-e^{-k(t-\Lambda_{0,max})}\right]\geq \frac{B}{k}\left[1-e^{-k\rho h}\right]=S^{\vartriangle}.
 \end{equation}
 Hence, from (\ref{ch1.sec4.lemma1.proof.eq6}) and (\ref{ch1.sec4.lemma1.proof.eq8}), it follows that for some suitable choice of $\rho>0$ sufficiently large, then
 \begin{equation}\label{ch1.sec4.lemma1.proof.eq9}
   S^{\vartriangle}>S^{*}_{1}, \forall t\geq \Lambda_{0,max}+ \rho h.
 \end{equation}

 For $t\geq \Lambda_{0,max}+ \rho h$, define
 \begin{eqnarray}
   V(t) &=& I(t) + \beta S^{*}_{1}\int_{t_{0}}^{h_{2}}\int_{t_{0}}^{h_{1}}f_{T_{2}}(u)f_{T_{1}}(s)e^{-\mu(s+u)}\int_{t-s}^{t}G(I(v-u))dvdsdu \nonumber\\
    &&+ \beta S^{*}_{1}\int_{t_{0}}^{h_{2}}\int_{t_{0}}^{h_{1}}f_{T_{2}}(u)f_{T_{1}}(s)e^{-\mu(s+u)}\int_{t-u}^{t}G(I(v))dvdsdu. \label{ch1.sec4.lemma1.proof.eq10}
 \end{eqnarray}
 It is easy to see from system (\ref{ch1.sec0.eq3})-(\ref{ch1.sec0.eq5}), and (\ref{ch1.sec4.lemma1.proof.eq10}) that differentiating $V(t)$ with respect to the system (\ref{ch1.sec0.eq3}) and (\ref{ch1.sec0.eq5}), leads to the following
 \begin{eqnarray}
   \dot{V}(t) &=&\beta \int_{t_{0}}^{h_{2}}\int_{t_{0}}^{h_{1}}f_{T_{2}}(u)f_{T_{1}}(s)e^{-\mu(s+u)}G(I(t-s-u))[S(t-u)-S^{*}_{1}]dsdu  \nonumber\\
   &&+\left[\beta S^{*}_{1}E(e^{-\mu(T_{1}+T_{2})})\frac{G(I(t))}{I(t)}-(\mu + d+\alpha)\right]I(t).\label{ch1.sec4.lemma1.proof.eq11}
 \end{eqnarray}
 For all $t\geq \Lambda_{0,max}+ \rho h +h>\Lambda_{0,max}+ \rho h +h_{2}$, it follows from (\ref{ch1.sec4.lemma1.eq1a}), (\ref{ch1.sec4.lemma1.proof.eq9}) and (\ref{ch1.sec0.eq3})- (\ref{ch1.sec0.eq5}) that
 \begin{eqnarray}
   \dot{V}(t) &\geq&\beta \int_{t_{0}}^{h_{2}}\int_{t_{0}}^{h_{1}}f_{T_{2}}(u)f_{T_{1}}(s)e^{-\mu(s+u)}G(I(t-s-u))[S^{\vartriangle}-S^{*}_{1}]dsdu \nonumber \\
   &&+\left[\beta S^{*}_{1}E(e^{-\mu(T_{1}+T_{2})})\frac{G(I^{*}_{1})}{I^{*}_{1}}-(\mu + d+\alpha)\right]I(t)\nonumber \\
   &=&\beta \int_{t_{0}}^{h_{2}}\int_{t_{0}}^{h_{1}}f_{T_{2}}(u)f_{T_{1}}(s)e^{-\mu(s+u)}G(I(t-s-u))[S^{\vartriangle}-S^{*}_{1}]dsdu.\label{ch1.sec4.lemma1.proof.eq11b}
 \end{eqnarray}
 Observe that the union of the subintervals  $\bigcup _{(s, u)\in [t_{0}, h_{1}]\times[t_{0}, h_{2}]}{[t_{0}-(s+u), t_{0}]}=[t_{0}-h, t_{0}]$, where $ h=h_{1}+h_{2}$. Denote the following
 \begin{equation}\label{ch1.sec4.lemma1.proof.eq11c}
   i_{min}=\min_{\theta\in [t_{0}-h, t_{0}], (s, u)\in [t_{0}, h_{1}]\times[t_{0}, h_{2}]} {I(\Lambda_{0,max}+ \rho h+ h+s+u+\theta)}.
 \end{equation}
 Note that (\ref{ch1.sec4.lemma1.proof.eq11c}) is equivalent to
 \begin{equation}\label{ch1.sec4.lemma1.proof.eq11a}
   i_{min}=\min_{\theta\in [t_{0}-h, t_{0}]} {I(\Lambda_{0,max}+ \rho h+ h+h+\theta)}.
 \end{equation}
 It is shown in the following that $I(t)\geq i_{min}, \forall t\geq \Lambda_{0,max}+ \rho h+h\geq \Lambda_{0,max}+ \rho h+u$, $\forall u\in [t_{0}, h_{2}]$.

 Suppose on the contrary there exists $\tau_{1}\geq 0$ such that  $I(t)\geq i_{min}$ for all
  $t\in [\Lambda_{0,max}+ \rho h +h, \Lambda_{0,max}+ \rho h +h+h+\tau_{1}]\supset [\Lambda_{0,max}+ \rho h +h, \Lambda_{0,max}+ \rho h +h+ s+u+\tau_{1}], \forall (s, u)\in [t_{0}, h_{1}]\times[t_{0}, h_{2}]$
 \begin{equation}\label{ch1.sec4.lemma1.proof.eq12}
   I(\Lambda_{0,max}+ \rho h +h+h+\tau_{1})=i_{min},\quad and\quad  \dot{I}(\Lambda_{0,max}+ \rho h +h+h+\tau_{1})\leq 0.
 \end{equation}
 For the value of $t=\Lambda_{0,max}+ \rho h +h+h+\tau_{1}$, it follows that $ S(t-u)>S^{\vartriangle}>S^{*}_{1}$, and $t-s-u\in [\Lambda_{0,max}+ \rho h +h, \Lambda_{0,max}+ \rho h +h+h+\tau_{1}]$, $ \forall  (s, u)\in [t_{0}, h_{1}]\times[t_{0}, h_{2}]$, and  it can be further seen from (\ref{ch1.sec0.eq3})-(\ref{ch1.sec0.eq5}), (\ref{ch1.sec4.lemma1.proof.eq9}) and (\ref{ch1.sec4.lemma1.eq1a}) that
 \begin{eqnarray}
   &&\dot{I}(t)\geq \beta E(e^{-\mu(T_{1}+T_{2})})G(i_{min})S^{\vartriangle}-(\mu+d+\alpha)i_{min}=\left[\beta E(e^{-\mu(T_{1}+T_{2})})\frac{G(i_{min})}{i_{min}}S^{\vartriangle}-(\mu+d+\alpha)\right]i_{min} \nonumber \\
    &&>\left[\beta E(e^{-\mu(T_{1}+T_{2})})\frac{G(I^{*}_{1})}{I^{*}_{1}}S^{*}_{1}-(\mu+d+\alpha)\right]i_{min}=0.\label{ch1.sec4.lemma1.proof.eq13}
 \end{eqnarray}
But (\ref{ch1.sec4.lemma1.proof.eq13}) contradicts (\ref{ch1.sec4.lemma1.proof.eq12}). Therefore, $I(t)\geq i_{min}, \forall t\geq \Lambda_{0,max}+ \rho h+h\geq \Lambda_{0,max}+ \rho h+u+s$, $\forall  (s, u)\in [t_{0}, h_{1}]\times[t_{0}, h_{2}]$.

 It follows further  from (\ref{ch1.sec4.lemma1.proof.eq11})-(\ref{ch1.sec4.lemma1.proof.eq11c}), and the Assumption~\ref{ch1.sec0.assum1} that for $\forall t\geq \Lambda_{0,max}+ \rho h+h+h\geq \Lambda_{0,max}+ \rho h+h+s+u$, $\forall (s,u)\in [t_{0}, h_{1}]\times[t_{0}, h_{2}]$.
 \begin{eqnarray}
   \dot{V}(t) &\geq&\beta \int_{t_{0}}^{h_{2}}\int_{t_{0}}^{h_{1}}f_{T_{2}}(u)f_{T_{1}}(s)e^{-\mu(s+u)}G(I(t-s-u))[S^{\vartriangle}-S^{*}_{1}]dsdu\nonumber\\
   &>& \beta E(e^{-\mu(T_{1}+T_{2})})G(i_{min})(S^{\vartriangle}-S^{*}_{1})>0.\label{ch1.sec4.lemma1.proof.eq14}
 \end{eqnarray}
 From (\ref{ch1.sec4.lemma1.proof.eq14}), it implies that $\limsup_{t\rightarrow\infty}{V(t)}=+\infty$.

  On the contrary, it can be seen from  (\ref{ch1.sec1.thm1a.eq0}) that  $\limsup_{t\rightarrow \infty} N(t)\leq S^{*}_{0}=\frac{B}{\mu}$, which implies that $\limsup_{t\rightarrow \infty} I(t)\leq S^{*}_{0}=\frac{B}{\mu}$. This further implies that for every $\epsilon>0$ infinitesimally small, there exists $\tau_{2}>0$ sufficiently large such that $I(t)\leq S^{*}_{0}+\varepsilon, \forall t\geq \tau_{2}$.  It follows that from Assumption~\ref{ch1.sec0.assum1} that
  \begin{equation}\label{ch1.sec4.lemma1.proof.eq15}
    G(I(t-s-u))\leq G(I(v-u))\leq G(I(t-u))\leq G(I(t))\leq G(S^{*}_{0}+\epsilon), \forall v\in [t-s,t], (s,u)\in [t_{0}, h_{1}]\times[t_{0}, h_{2}].
  \end{equation}
  From (\ref{ch1.sec4.lemma1.proof.eq15}), it follows that
  \begin{equation}\label{ch1.sec4.lemma1.proof.eq16}
  \limsup_{t\rightarrow\infty}{G(I(t-s-u))}\leq \limsup_{t\rightarrow\infty}{G(I(t))}\leq G(S^{*}_{0}).
  \end{equation}
  It is easy to see from (\ref{ch1.sec4.lemma1.proof.eq10}) and (\ref{ch1.sec4.lemma1.proof.eq16}) that
  \begin{equation}\label{ch1.sec4.lemma1.proof.eq17}
    \limsup_{t\rightarrow\infty}{V(t)}\leq S^{*}_{0}+\beta S^{*}_{1}G(S^{*}_{0})E\left((T_{1}+T_{2})e^{-\mu(T_{1}+T_{2})}\right)<\infty.
  \end{equation}
  Therefore, it is impossible that $I(t)\leq q I^{*}_{1}$ for sufficiently large  $t\geq t_{0}$, where $q\in(0, 1)$.

  Hence, the following are possible, $(Case(i.))$ $I(t)\geq qI^{*}_{1}$ for all $t$ sufficiently large, and $(Case(ii.))$ $I(t)$ oscillates about $qI^{*}_{1}$  for sufficiently large $t$. Obviously, we need show only $Case(ii.)$. Suppose $t_{1}$ and $t_{2}$ are are sufficiently large values such that
  \begin{equation}\label{ch1.sec4.lemma1.proof.eq18}
    I(t_{1})=I(t_{2})= qI^{*}_{1},\quad and\quad I(t)<qI^{*}_{1}, \forall (t_{1}, t_{2}).
  \end{equation}
  If for all $(s,u)\in [t_{0}, h_{1}]\times[t_{0}, h_{2}]$,  $t_{2}-t_{1}\leq \rho h+h$, where $ h=h_{1}+h_{2}$, observe that $[t_{1}, t_{1}+\rho h+s+u]\subseteq [t_{1}, t_{1}+\rho h+h]$, and  it is easy to see from (\ref{ch1.sec0.eq3}) by integration that
  \begin{equation}\label{ch1.sec4.lemma1.proof.eq19}
    I(t)\geq I(t_{1})e^{-(\mu+d+\alpha)(t-t_{1})}\geq qI^{*}_{1} e^{-(\mu+d+\alpha)(\rho+1)h}\equiv v_{2}.
  \end{equation}
 If for all $(s,u)\in [t_{0}, h_{1}]\times[t_{0}, h_{2}]$,  $t_{2}-t_{1}>\rho h+h\geq \rho h+ s+u$, then it can be seen easily that $I(t)\geq v_{2}$, for all $t\in [t_{1}, t_{1}+\rho h+s+u]\subseteq  [t_{1}, t_{1}+\rho h+h]$.

 Now, for each $t\in (\rho h+h, t_{2})\supseteq (\rho h+s+u, t_{2})$, $\forall (s,u)\in [t_{0}, h_{1}]\times[t_{0}, h_{2}]$,
  one can also claim that $I(t)\geq v_{2}$.  Indeed, as similarly shown above, suppose  on the  contrary for all $(s,u)\in [t_{0}, h_{1}]\times[t_{0}, h_{2}]$, $\exists T^{*}>0$ such that $I(t)\geq v_{2}$, $\forall t\in [t_{1}, t_{1}+\rho h+h+T^{*}]\supseteq [t_{1}, t_{1}+\rho h+s+u+T^{*}]$
 \begin{equation}\label{ch1.sec4.lemma1.proof.eq20}
   I(t_{1}+\rho h+h+T^{*})=v_{2},\quad but \quad \dot{I}(t_{1}+\rho h+h+T^{*})\leq 0.
 \end{equation}
 It follows from (\ref{ch1.sec0.eq3})-(\ref{ch1.sec0.eq5}) and (\ref{ch1.sec4.lemma1.eq1a}) that for the value of $t=t_{1}+\rho h+h+T^{*}$, 
 \begin{eqnarray}
   &&I(t) \geq \beta E(e^{-\mu(T_{1}+T_{2})})G(v_{2})S^{\vartriangle}-(\mu+d+\alpha)v_{2}>\left[\beta E(e^{-\mu(T_{1}+T_{2})})\frac{G(v_{2})}{v_{2}}S^{*}_{1}-(\mu+d+\alpha)\right]v_{2} \nonumber \\
    &&\geq\left[\beta E(e^{-\mu(T_{1}+T_{2})})\frac{G(I^{*}_{1})}{I^{*}_{1}}S^{*}_{1}-(\mu+d+\alpha)\right]v_{2}=0.\label{ch1.sec4.lemma1.proof.eq21}
 \end{eqnarray}
  Observe that (\ref{ch1.sec4.lemma1.proof.eq21}) contradicts (\ref{ch1.sec4.lemma1.proof.eq20}). Therefore, $I(t)\geq v_{2}$,  for $t\in [t_{1}, t_{2}]$.
  And since  $[t_{1}, t_{2}]$ is arbitrary, it implies that $I(t)\geq v_{2}$ for all sufficiently large $t$.  Therefore (\ref{ch1.sec4.lemma1.eq1}) is satisfied.
%
 \begin{thm}\label{ch1.sec4.thm1}
 If the conditions of Lemma~ \ref{ch1.sec4.lemma1} are satisfied, then the system (\ref{ch1.sec0.eq3})-(\ref{ch1.sec0.eq5}) is strongly permanent for any total delay time $h=h_{1}+h_{2}$  according to Definition~\ref{ch1.sec0.eq13b.def1}[1]. 
 \end{thm}
  \begin{rem}\label{ch1.sec4.rem1}
  \item[1.] It can be seen from Lemma~ \ref{ch1.sec4.lemma1} (\ref{ch1.sec4.lemma1.eq1}) that when $\beta=0$, then $v_{1}=\frac{B}{\mu}$. That is, when disease transmission stops, then asymptotically, the smallest total susceptible that remains are new births over the average lifespan $\frac{1}{\mu}$ of the population, equivalent to the DFE $S^{*}_{0}=\frac{B}{\mu}\equiv 1$. Also, as  $\beta\rightarrow\infty$, then the total susceptible that remains $v_{1}\rightarrow 0^{+}$. That is, as disease transmission rises, even the new births are either infected, or die from natural or disease related causes over time.
    \item[2.] From (\ref{ch1.sec4.lemma1.eq1}), observe that $e^{-(\mu+d+\alpha)(\rho+1)h}$ is the survival probability  from natural death ($\mu$), disease mortality ($d$), and from infectiousness ($\alpha$), over the total life cycle of the parasite $h$. Thus, the smallest infectious state  that remains asymptotically $v_{2}\equiv qI^{*}_{1}e^{-(\mu+d+\alpha)(\rho+1)h}$ is a fraction $q\in(0,1)$ of the endemic equilibrium $I^{*}_{1}$ that survives from death and disease over life cycle $h$.

        Since $\frac{1}{(\mu+d+\alpha)}$ is the effective average lifespan of an individual who survives the disease until recovery at rate $\alpha$, it follows from (\ref{ch1.sec4.lemma1.eq1}) that as $(\mu+d+\alpha)\rightarrow 0^{+}$ and consequently $\frac{1}{(\mu+d+\alpha)}>>1$, then $e^{-(\mu+d+\alpha)(\rho+1)h}\rightarrow 1^{-}$. Moreover, from (\ref{ch1.sec4.lemma1.eq1}), $v_{2}\rightarrow qI^{*}_{1}$. That is, when malaria is actively transmitted $\beta>0$, so that more susceptible individuals become infected, but the effective average lifespan is still high because for example, malaria is  treated, or healthier lifestyles are encouraged, and less people die from the disease  $0\leq d<<1$, and from natural causes $0\leq\mu<<1$, then the total infectious state that remains over time $v_{2}$ is a fraction $q>>qe^{-(\mu+d+\alpha)(\rho+1)h}$ of all infected at steady state $I^{*}_{1}$, i.e. more  infectious remains over time when malaria is treated effectively, or healthier living standards are encouraged.
        \item[3.] The question of what conditions the population ever gets extinct in time is answered from [1.]~\&~[2.] above. Since as $\beta\rightarrow +\infty$, and $(\mu+d+\alpha)\rightarrow \infty$,  then $v_{1}\rightarrow 0^{+}$, and $v_{2}\rightarrow 0^{+}$, respectively, from (\ref{ch1.sec4.lemma1.eq1}). Thus,  extinction is ever possible in time, whenever disease transmission rate is high, and the response to malaria treatment or living standards are very poor. 
        %
 \end{rem}

 \section{Example: Application to P. vivax malaria}\label{ch1.sec4}
   In this section, the extinction and persistence results are exhibited for the \textit{P.vivax malaria} example in Wanduku \cite{wanduku-extinct}. This is accomplished by examining the  trajectories of the decoupled system (\ref{ch1.sec0.eq3}) and (\ref{ch1.sec0.eq5}) relative to the zero and endemic equilibria. To conserve space, we recall the dimensionless parameters in [Table 1,\cite{wanduku-extinct}, page 3793] given in Table~\ref{ch1.sec4.table1}, and the reader is referred to \cite{wanduku-extinct} for detailed description of the \textit{P.vivax malaria} scenario.

   The dimensional estimates for the parameters of the malaria model given in [(a.)-(e.), \cite{wanduku-extinct}, page 3792] are applied to (\ref{ch1.sec0.eq6.eq20}) to find the dimensionless parameters for the model (\ref{ch1.sec0.eq3})-(\ref{ch1.sec0.eq6}) given in Table~\ref{ch1.sec4.table1}.

    \begin{table}[h]
  \centering
  \caption{A list of dimensionless values  for the system parameters for Example 1. }\label{ch1.sec4.table1}
  \begin{tabular}{l l l}
  Disease transmission rate&$\beta$& Subsection~\ref{ch1.sec4.subsec1} ($0.02146383$), Subsection~\ref{ch1.sec4.subsec2} ($0.2146383$)\\\hline
  Constant Birth rate&$B$&$ 8.476678e-06$\\\hline
  Recovery rate& $\alpha$& $0.08571429$\\\hline
  Disease death rate& $d$& 0.0001761252\\\hline
  Natural death rate& $\mu$, $\mu_{v}$& $ 8.476678e-06, 42.85714$\\\hline
  Incubation delay in vector& $T_{1}$& 0.105 \\\hline
  Incubation delay in host& $T_{2}$& 0.175\\\hline
  Immunity delay time& $T_{3}$& 2.129167\\\hline
  \end{tabular}
\end{table}

    Moreover, the  Euler approximation scheme is used to generate trajectories for the different states $S(t), E(t), I(t), R(t)$ over the time interval $[0,1000]$ days. The special nonlinear incidence  functions $G(I)=\frac{a_{1}I}{1+I}, a_{1}=0.05$ in \cite{gumel} is utilized. Furthermore, the following initial fractions of susceptible, exposed, infectious and removed individuals in the initial population size $\hat{N}(t_{0})=65000$ are used:
\begin{eqnarray}
&&S(t)= 10/23\approx 28261/65000, E(t)= 5/23\approx 14131/65000, I(t)= 6/23\approx 16957/65000,\nonumber\\
&&R(t)= 2/23\approx 5653/65000,\forall t\in [-T,0], T=\max(T_{1}+T_{2}, T_{3})=2.129167.\label{ch1.sec4.eq1}
\end{eqnarray}

Recall Section~\ref{ch1.sec1} asserts that the endemic equilibrium $E_{1}$ exists, whenever the BRN $R^{*}_{0}>1$, where $R^{*}_{0}$ is defined in (\ref{ch1.sec2.lemma2a.corrolary1.eq4}). Thus, it follows that when $R^{*}_{0}>1$, the endemic equilibrium $E_{1}=(S^{*}_{1}, E^{*}_{1},I^{*}_{1}, R^{*}_{1})$ satisfies the following system
\begin{eqnarray}
&&B-\beta Se^{-\mu_{v} T_{1}}G(I)-\mu S+\alpha I e^{-\mu T_{3}}=0, \beta Se^{-\mu_{v} T_{1}}G(I)-\mu E -\beta Se^{- (\mu_{v}T_{1}+\mu T_{2})}G(I)=0,\nonumber\\
&&\beta Se^{- (\mu_{v}T_{1}+\mu T_{2})}G(I)-(\mu+d+\alpha)I=0,\alpha I-\mu R-\alpha I e^{-\mu T_{3}}=0.\label{ch1.sec4.eq1}
\end{eqnarray}

For the given set of dimensionless parameter estimates in Table~\ref{ch1.sec4.table1}, the DFE is $E_{0}=(S^{*}_{0}, 0,0)=(1,0,0)$. Also, the endemic equilibrium is given as $E_{1}=(S^{*}_{1}, E^{*}_{1},I^{*}_{1})=(0.002323845,0.00068247,0.04540019)$.

 \subsection{Example for extinction of disease}\label{ch1.sec4.subsec1}
 For the given set of dimensionless parameter estimates in Table~\ref{ch1.sec4.table1}, where $\beta=0.02146383$, from (\ref{ch1.sec2.lemma2a.corrolary1.eq4}) the BRN is $\hat{R}^{*}_{0}=0.2498732<1$. Therefore, $E_{0}$ is stable, and the endemic equilibrium $E_{1}=(S^{*}_{1}, E^{*}_{1},I^{*}_{1})$ fails to exist.
\begin{figure}[H]
\begin{center}
\includegraphics[scale=0.2]
{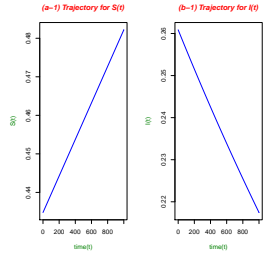}
\caption{(a-1), and (b-1),  show the trajectories of the  states $(S,I)$, respectively, over sufficiently long time $t\in [0,1000]$, whenever the intensity of the incidence of malaria is $a=0.05$.
  The BRN in (\ref{ch1.sec2.lemma2a.corrolary1.eq4}) in this case is $R^{*}_{0}= 0.2498732<1$, the estimate of the LE, or rate of extinction of the disease in (\ref{ch1.sec2b.thm1.eq1.proof.eq2.eq1}) is $\lambda= 0.06443506>0$.
   }\label{ch1.sec4.subsec1.fig2}
\end{center}
\end{figure}

  Figure~\ref{ch1.sec4.subsec1.fig2} verifies the results about the extinction of the $I(t)$ state over time in Theorem~\ref{ch1.sec2b.thm1}, and the persistence of the $S(t)$ state over timein Theorem~\ref{ch1.sec2b.thm2}. Indeed, it is observed that for the given parameter values in Table~\ref{ch1.sec4.table1}, and the initial conditions  in (\ref{ch1.sec4.eq1}),  the BRN in (\ref{ch1.sec2.lemma2a.corrolary1.eq4}) in this scenario is $R^{*}_{0}= 0.2498732<1$. Therefore, the condition of Theorem~\ref{ch1.sec2b.thm1}(a.) and Theorem~\ref{ch1.sec2b.thm2} are satisfied, and from (\ref{ch1.sec2b.thm1.eq1.proof.eq2.eq1}), the estimate of the rate of extinction of the malaria population $I(t)$ is $\lambda= 0.06443506>0$. That is,
\begin{equation}\label{ch1.sec4.subsec1.eq1}
  \limsup_{t\rightarrow \infty}{\frac{1}{t}\log{(I(t))}}\leq -\lambda = -0.06443506.
\end{equation}
The Figure~\ref{ch1.sec4.subsec1.fig2}(b-1) confirms that over sufficiently large time, when $\lambda>0$, then the infectious state approaches zero, that is, $\lim_{t\rightarrow \infty}I(t)=0$. Furthermore, the BRN  $R^{*}_{0}= 0.2498732<1$,  signifies that the disease is getting eradicated from the population over time. This is confirmed by Figure~\ref{ch1.sec4.subsec1.fig2}(a-1), where $S(t)$ appears to be rising over time, and approaching the DFE state $S^{*}_{0}=\frac{B}{\mu}=1$, that is, $\lim_{t\rightarrow \infty}S(t)=1$.
%
\subsection{Persistence of malaria}\label{ch1.sec4.subsec2}
For the given set of dimensionless parameter estimates in Table~\ref{ch1.sec4.table1}, when $\beta=7.941616$,  from (\ref{ch1.sec2.lemma2a.corrolary1.eq4}) the BRN becomes $\hat{R}^{*}_{0}=92.45307>1$. Therefore, the DFE $E_{0}=(S^{*}_{0}, 0,0)=(1,0,0)$ becomes unstable, and   the endemic equilibrium exists, and given as $E_{1}=(S^{*}_{1}, E^{*}_{1},I^{*}_{1})=(6.281296e-06,0.0006840553,0.04550565)$.

\begin{figure}[H]
\begin{center}
\includegraphics[scale=0.2]
{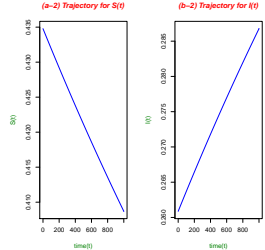}
\caption{(a-2), and (b-2),  show the trajectories of the  states $(S,I)$, respectively, over sufficiently long time $t\in [0,1000]$, whenever the intensity of the incidence of malaria is $a=0.05$.
  The BRN in (\ref{ch1.sec2.lemma2a.corrolary1.eq4}) is $R^{*}_{0}= 92.45307<1$. $E(e^{- (\mu_{v}T_{1}+\mu T_{2})})=0.01110898>0.0108163=\frac{1}{R^{*}_{0}}$. Hence, Theorem~\ref{ch1.sec2b.thm1}  and Theorem~\ref{ch1.sec2b.thm2} fail.  The endemic equilibrium $E_{1}=(S^{*}_{1}, E^{*}_{1},I^{*}_{1})$ exists and Theorem~\ref{ch1.sec4.thm1} holds from (\ref{ch1.sec4.subsec2.eq1}).
   }\label{ch1.sec4.subsec1.fig3}
\end{center}
\end{figure}

It can be shown from (\ref{ch1.sec4.lemma1.eq1}) that for some suitable choice of $q\in (0,1)$ and $\rho>0$, for $t\in [0,1000]$,
 \begin{eqnarray}
    \liminf_{t\rightarrow \infty}{S(t)}=0.4086943>> v_{1}\equiv \frac{B}{\mu+\beta G(S_{0})}=4.269316e-05,\nonumber\\
    \liminf_{t\rightarrow \infty}{I(t)}=0.2360496>> v_{2}\equiv qI^{*}_{1}e^{-(\mu+d+\alpha)(\rho+1)h}=0.04550565qe^{-0.02405169(1+\rho)}.\label{ch1.sec4.subsec2.eq1}
  \end{eqnarray}
Hence, from Theorem~\ref{ch1.sec4.thm1} there is a significant number of infectious people present over time $[0,1000]$, and as a result malaria persists in the population over time. These facts are further illustrated by Figure~\ref{ch1.sec4.subsec1.fig3} over $[0,1000]$.

 \section{conclusion}
  The vector-human population dynamic models are derived. The models have a general nonlinear incidence rate. The extinction and persistence of the vector-borne disease in the SEIRS epidemic models are studied. Numerical simulation results are given to confirm the results.
\paragraph{Acknowledgment.}
Thanks to the Editor and reviewers for the thorough and constructive feedback.

\bibliographystyle{NAplain}
\bibliography{12-03-2019-references-sample}
\end{document}